\documentclass[11pt,a4paper]{amsart}
\usepackage{amsfonts,amssymb,bbm,enumerate,mathrsfs}

\def\C{\Bbb{C}}

\def\N{\Bbb{N}}

\def\R{\Bbb{R}}

\def\di{\partial}

\def\bl{\langle}
\def\br{\rangle}

\def\liml{\lim\limits}
\def\suml{\sum\limits}\def\intl{\int\limits}
\def\maxl{\max\limits}
\def\oplusl{\mathop\oplus\limits}

\def\cupl{\mathop\cup\limits}

\def\empty{\varnothing}

\newcommand{\quot}[2]{{\left.\raisebox{1.6ex}{\footnotesize$#1$}  \!\!\!\!{\ensuremath\diagup}\!\!\raisebox{-1ex}{\footnotesize$#2$}\right.}}
\newcommand{\quots}[2]{{\footnotesize\left.\raisebox{0.4ex}{$#1$}\! / \!\raisebox{-0.4ex}{$#2$}\right.}}

\def\tA{\tilde{A}}\def\tf{\tilde{f}}
\def\tF{\tilde{F}}\def\tg{\tilde{g}}\def\th{\tilde{h}}\def\tI{\tilde{I}}\def\tk{{\tilde{k}}}
\def\tR{{\tilde{R}}}\def\ttau{\tilde{\tau}}

\def\hA{\widehat{A}}\def\hF{\widehat{F}}\def\hf{\widehat{f}}\def\hg{{\widehat{g}}}\def\hI{{\widehat{I}}}
\def\hM{{\widehat{M}}}\def\hR{{\widehat{R}}}
\def\hy{\widehat{y}}

\def\al{\alpha}
\def\ep{\epsilon}
\def\om{{\omega}}

\def\ca{\frak a}
\def\cb{\frak b}
\def\cc{\frak c}

\def\cm{{\frak m}}

\def\cq{{\frak q}}
\def\cU{\mathcal U}
\def\cV{\mathcal V}

\def\uk{{\underline{k}}}
\def\ul{{\underline{l}}}\def\um{{\underline{m}}}

\def\ux{{\underline{x}}}\def\uz{{\underline{z}}}

\def\one{{1\hspace{-0.1cm}\rm I}}

\newcommand{\ber}{\begin{array}{l}}\newcommand{\eer}{\end{array}}
\newcommand{\bpm}{\begin{pmatrix}}\newcommand{\epm}{\end{pmatrix}}
\newcommand{\bbm}{\begin{bmatrix}}\newcommand{\ebm}{\end{bmatrix}}
\newcommand{\bM}{\begin{matrix}}\newcommand{\eM}{\end{matrix}}
\newcommand{\bee}{\begin{enumerate}}\newcommand{\eee}{\end{enumerate}}
\newcommand{\bei}{\begin{itemize}}\newcommand{\eei}{\end{itemize}}

\def\wrt{with respect to }
\def\sset{\subset}\def\sseteq{\subseteq}\def\ssetneq{\subsetneq}\def\smin{\setminus}

\def\Mat{\operatorname{Mat}_{m\times n}(R)}

\newtheorem{Lemma}{Lemma}[section]\newcommand{\bel}{\begin{Lemma}}\newcommand{\eel}{\end{Lemma}}
\newtheorem{Theorem}[Lemma]{Theorem}\newcommand{\bthe}{\begin{Theorem}}\newcommand{\ethe}{\end{Theorem}}
\newtheorem{Proposition}[Lemma]{Proposition}\newcommand{\bprop}{\begin{Proposition}}\newcommand{\eprop}{\end{Proposition}}
\newtheorem{Corollary}[Lemma]{Corollary}\newcommand{\bcor}{\begin{Corollary}}\newcommand{\ecor}{\end{Corollary}}

\newtheorem{Definition}[Lemma]{Definition}
\newcommand{\bed}{\begin{Definition}}
\newcommand{\eed}{\end{Definition}}
\newtheorem{Definition-Proposition}[Lemma]{Definition-Proposition}

\def\bpr{~\\{\em Proof.\ \ }}
\newcommand{\epr}{{\hfill\ensuremath\blacksquare}}
\newtheorem{Remark}[Lemma]{Remark}
\newcommand{\beR}{\begin{Remark}\rm}
\newcommand{\eeR}{\end{Remark}}
\newtheorem{Example}[Lemma]{Example}
\newcommand{\bex}{\begin{Example}\rm}
\newcommand{\eex}{\end{Example}}
\newtheorem{Problem}[Lemma]{Problem}
\newcommand{\bprob}{\begin{Problem}\rm}
\newcommand{\eprob}{\end{Problem}}

\newcommand{\bet}{\begin{tabular}{cccccccc}}\newcommand{\eet}{\end{tabular}}
\newcommand{\beq}{\begin{equation}}\newcommand{\eeq}{\end{equation}}

\newcommand{\bin}[2]{\binom{#1}{#2}}
\newcommand\isom{\xrightarrow{\,\smash{\raisebox{-0.65ex}{\ensuremath{\scriptstyle\sim}}}\,}}

\title[]{S\MakeLowercase{urjectivity of the completion map for rings of $\MakeUppercase{C}^\infty$-functions.}
\\(W\MakeLowercase{hitney extension theorem for general filtrations)} }
\author[]{G\MakeLowercase{enrich} B\MakeLowercase{elitskii and}
D\MakeLowercase{mitry} K\MakeLowercase{erner}}
\address{Department of Mathematics, Ben Gurion University of the Negev, P.O.B. 653, Be'er Sheva 84105, Israel.}
\email{genrich@math.bgu.ac.il}
\email{dmitry.kerner@gmail.com}
\date{\today\ \  filename: \jobname.tex}
\thanks{D.Kerner was partially supported by Israel Science Foundation (grant No.  1910/18)}


\begin{document}
\begin{abstract}
The classical lemma of Borel reads: any power series with real coefficients is
 the Taylor series of a smooth function.
 Algebraically this means the surjectivity of the completion map at a point,
 $C^\infty(\R^n)\twoheadrightarrow\R[[\ux]]$.
  Similarly, Whitney extension theorem implies the surjectivity of
  the completion at closed subsets of $\R^n$.

 For various applications one needs the surjectivity of completion for general $C^\infty$-rings and general filtrations.
   We establish the necessary and sufficient conditions for this surjectivity.

 Moreover, we prove: any element of the completion admits a $C^\infty$-representative that is real-analytic outside of the locus of completion,
 has any prescribed vanishing rate ``at infinity", and the
  prescribed positivity behaviour at the finite part.
  Alternatively, one can impose on the smooth representative a  set of (compatible) linear conditions.
\end{abstract}
\maketitle
\setcounter{secnumdepth}{6} \setcounter{tocdepth}{1} 

\section{Introduction}
\subsection{}
The classical lemma of  Borel\footnote{Though published in \cite{Borel}, it was partially known to Peano, see \cite{Besenyei}.} reads:
 any sequence of real numbers is realizable as the sequence of partial derivatives (at the origin) of a smooth function,
 moreover this function can be assumed analytic off the origin.
 Algebraically this means the surjectivity of the completion map $C^\infty(\R^n)\twoheadrightarrow\R[[\ux]]$, for the filtration $(\ux)^\bullet$. We
  get two exact sequences:
\beq
\bM 0\to (\ux)^\infty\to& C^\infty(\R^n) &\to &\R[[\ux]]&\to 0
\\&
\cup
&
&||&\\
& C^\infty(\R^n)\cap C^\om(\R^n\smin \{o\})& \to &\R[[\ux]]&\to 0
\eM
\eeq
This surjectivity goes in notable difference to the completion maps of traditional (Noetherian) rings of Commutative Algebra/Algebraic Geometry.

 More generally, Whitney extension theorem
  gives the necessary and sufficient conditions to extend a function with prescribed ``derivatives" on a(ny) closed
   set $Z\sset \cU\sseteq \R^n$ to a  function smooth on $\cU$,  \S\ref{Sec.Whitney.Extension.Theorem}.
  (Moreover, the function can be taken real analytic on $\cU\smin Z$.)
 Algebraically this implies the two exact sequences of completion:
\beq\label{Eq.Whitney.Thm.Completion}
\bM 0\to I(Z)^{\bl \infty\br}\to& C^\infty(\cU) &\to& \widehat{C^\infty(\cU)}^{(Z)}&\to 0
\\& \cup
&& || \\
&  C^\infty(\R^n)\cap C^\om(\cU\smin Z)& \to &\widehat{C^\infty(\cU)}^{(Z)}&\to 0
\eM
\eeq
(Here $I(Z)^{\bl \infty\br}$ denotes the ideal of functions flat on $Z$.
  The completion is taken \wrt the differential powers of ideals, we denote this filtration by $I(Z)^{\bl \bullet\br}$, see \S\ref{Sec.Notations.Conventions}.ii.)

\

These results are of everyday use in Analysis/Differential Geometry. In Algebraic Geometry over $C^\infty$-rings, \cite{Joyce}, \cite{Moerdijk-Reyes},
  and in $C^\infty$-Singularity Theory
  one uses more general $C^\infty$-rings and more general filtrations, $I_\bullet$, not necessarily by powers of ideals.
    The surjectivity of completion does not always hold in this case,
   and the (necessary/sufficient) conditions are not implied by Whitney extension theorem.
  The data of derivatives of a function is essentially different from the data of an element of $\widehat{C^\infty(\cU)}\ \!^{(I_\bullet)}$.
Though this surjectivity question is most natural and basic,  it was
not   addressed before. (Perhaps because the methods needed involve both analysis and commutative algebra.)
 In \cite{Bel.Boi.Ker} we have obtained the first sufficient condition for this surjectivity,
  but it was far from being necessary.

\

Our paper extends the classical results in several directions:
\bei
\item We give the {\em necessary and sufficient} conditions for the surjectivity of the completion map for  general filtrations.
 Moreover, for an element of the completion $\hf\in \hR^{(I_\bullet)}$, we establish the preimage, $f\in R$, which is real analytic (in the maximal allowed region),
 decays with any prescribed rate at the boundary, $\di \cU$, and is positive off a small neighbourhood of $Z$.

Without requiring the analyticity on $\cU\smin Z$ we can ensure that $f$ satisfies any prescribed system of $\R$-linear conditions
 that are compatible in a precise sense. This extends the classical moment problem.
\item
 We establish this surjectivity for several types of $C^\infty$-rings, e.g. the quotient rings, $\quots{C^\infty(\cU)}{J}$,
   and localizations, $C^\infty(\cU)[S^{-1}]$. In particular,
    this addresses  the functions on manifolds/$C^\infty$-schemes, and the $C^\infty$-germs (at points or along closed sets).
\eei

 The sufficient condition (in \S4) has a different form   than the necessary condition (in \S3), but is trivially {\em weaker} than the necessary condition.
  Therefore both conditions are necessary and sufficient (and equivalent).
The necessary/sufficient conditions are stated algebraically (via the equivalence of filtrations by ideals).
 An interesting question is whether/how to formulate (and prove) these conditions  in terms of the classical analysis.

 We remark that the surjectivity conditions have nothing to do with $I_\bullet$
 being finitely/sub-analytically/$o$-minimally generated. See example \ref{Ex.Ideal.Polynomially.Generated}
   for polynomially generated $I_\bullet$ with non-surjective completion.
 It is impressive that this ``algebraization" of Whitney extension theorem is manageable in full generality,
 in view of the notoriously complicated/pathological behavior of (filtrations by) ideals in $C^\infty$-rings.

Finally, the additional properties of the preimage, $f\in R$, are valuable in various applications.

\

Due to the lack of space we postpone other results/applications (e.g. inversion of Borel-Taylor map for $I_\bullet$ and  the Tougeron-type approximations)
 to the next paper.

\

This surjectivity is the necessary ``Step.0" for various problems in Singularity Theory and Commutative Algebra.
 The standard approach is:   pass to the completion, $R\to \hR$,   resolve the problem over $\hR$, then   pullback the results to $R$.
 The later step is done via the Artin-Tougeron approximation, which in the $C^\infty$-case begins with ``Take a $C^\infty$-representative of $\hf\in \hR$."
 (See \cite{B.K.IFT}, \cite{Bel.Boi.Ker} for the $C^\infty$-approximation results and further references.)
 We were particularly motivated by the study of determinacy/algebraizability/deformations/local topology of $C^\infty$-maps/schemes
  with non-isolated singularities, \cite{B.K.motor}, \cite{BGK}, \cite{B.K.fin.det.2}.
   As the filtrations there are not by powers of ideals, the classical Whitney does not help, see e.g. \S\ref{Sec.Examples.surjectivity.of.completions}.vi.

\subsection{The structure of the paper}
 \bee[\S 1]\setcounter{enumi}{1}
\item is preparatory. In \S2.1 we fix the notations and conventions.
\bee[\S 2.1]\addtocounter{enumii}{1}
\item  
presents the classical Whitney extension theorem
 as the surjectivity of completion for the special filtration, by differential powers of ideals, $I(Z)^{\bl\bullet\br}$.
\item 
 recalls the analytic approximation of Whitney.
\item 
we prove the filtered (strengthened) version of Whitney theorem on zero sets.
 For any filtration $I_\bullet$ of $C^\infty(\cU)$, the total zero set $Z:=\overline{\cup  {V(I_j)}}\sset \cU$ is defined by one element of $I_\infty:=\cap I_j$.
  Moreover, this element can be chosen real-analytic on $\cU\smin Z$.
  This goes in notable difference to the traditional Noetherian rings of Algebraic Geometry, and can be considered as   a $C^\infty$-Nullstellensatz. It  is used repeatedly   later.

\item
 we show the persistence of surjectivity under the change of rings, for taking quotients,
 $R\to \quots{R}{J}$, and localizations,  $R\to R[S^{-1}]$. This reduces the surjectivity problem  to the
  very particular ring  $C^\infty(\cU)$, for an open $\cU\sseteq \R^n$, that can be chosen as a ball.

In particular, the surjectivity criteria for $C^\infty(\cU)$ apply, e.g. to the rings of $C^\infty$-germs (along closed sets).
\item  the surjectivity   for  $C^\infty$-rings on manifolds/schemes
 is reduced to the surjectivity  for
 $\quots{C^\infty(\cU)}{J}$. 
\eee

\item gives the necessary condition for the surjectivity, theorem \ref{Thm.Completion.Nec.Condition}.
 The proof goes by tracing the zero loci of the ideals $I_\bullet$,
 and the vanishing orders of $I_\bullet$ along the zero loci (\wrt the maximal ideals).
  Then we show the pointwise stabilization of the filtration $\quots{I_\bullet}{I_\infty}\sset \quots{R}{I_\infty}$.

A peculiar consequence is: if the $I_\bullet$-completion of $C^\infty(\cU)$ is surjective then it is non-injective, $I_\infty:=\cap I_j\neq(0)$.

 \item gives the sufficient condition for the surjectivity, theorem \ref{Thm.Completion.Suff.Condition}.
  Moreover,   we show that   any element of $\widehat{C^\infty(\cU)}^{\ (I_\bullet)}$ admits
   a representative $f\in C^\infty(\cU)$  that is real-analytic on $\cU\smin Z$, and positive outside of a small neighbourhood of $Z$.
  The positivity of $f$ cannot be strengthened to the positivity on the whole $\cU\smin Z$, see remark \ref{Re.Positivity.Cannot.Strengthen}.

  The proof is an explicit construction and (unsurprisingly) uses    cutoff functions with  derivatives of controlled growth.

\item contains    the simplest examples  and applications.
 In \S\ref{Sec.Examples.surjectivity.of.completions}  we derive  various classical
statements, e.g. multi-Borel lemma, Borel lemma in families, Borel
lemma for flat functions. In \S\ref{Sec.Examples.surjectivity.of.completions}.vi we consider a simple non-isolated
 hypersurface singularity. The corresponding filtration is not by powers of ideals, thus already in this
  case the needed surjectivity cannot be obtained with the classical Whitney extension theorem.

As   immediate applications of the surjectivity we address
 the ``Inverse Artin-Tougeron problem" in \S\ref{Sec.Inverse.Artin-Tougeron}, the contraction of ideals in \S\ref{Sec.Contraction.of.Ideals},
 the completion of $C^\infty$-modules in  \S\ref{Sec.Surjectivity.Completion.of.Modules}, and lifting the $\hR$-modules to $R$-modules in
  \S\ref{Sec.Lifting.of.Modules}.

\item        strengthens   the surjectivity further, by ensuring special properties of $f$.
\bee[\S6.1:]
\item
  a representative $f\in C^\infty(\cU)\cap C^\om(\cU\smin Z)$
 can be chosen with     any (prescribed) rate of decay on $\di\cU$. This strengthens/extends the Borel-type results
  for spaces of Schwartz functions. 
\item
   a representative $f\in C^\infty(\cU)$ can be chosen to satisfy
 any prescribed set of linear conditions   that are compatible in a precise sense.
  This extends numerous versions of the classical moment problem.
\eee
\eee

\subsection{Acknowledgement}
We thank M.Sodin for the highly useful reference to \cite{Hormander},
 A. Fernandez-Boix for the participation in the initial stage of this work, and  A. Kiro, E. Shustin for important comments.

\section{Preparations}\label{Sec.Preparations}

\subsection{Notations and conventions}\label{Sec.Notations.Conventions} Take an open subset $\cU\sseteq \R^n$ (possibly topologically
 non-trivial, e.g. non-contractible),
  an ideal $\{0\}\sseteq J\sset C^\infty(\cU)$, and the quotient ring $R=\quots{C^\infty(\cU)}{J}$.
\bee[\bf i.]
\item  For any ideal $I\sset  R$ we take  its (reduced) set of zeros,  $V(I)\sseteq V(J)\sseteq \cU$. For any  subset $Z\sset V(J)\sseteq\cU$ denote
 by $I(Z)\sset R$ the ideal of functions vanishing on $Z$. Thus $I(V(I))\supseteq I$.

 For any $\ep>0$ define the $\ep$-neighbourhood, $\cU_\ep(Z):=\{x|\ dist(x,Z)<\ep\}\sset \cU$.
\item
Denote the maximal ideal of a point $x\in \cU$ by $\cm_x\sset R$.

For any subset $Z\sset V(J)\sseteq \cU$ define the $j$'th differential power of $I(Z)$
 as the set of functions that vanish on $Z$ up to $j$'th order (see e.g. \cite{Symbolic.Differential.Powers})
\beq
I(Z)^{\bl j\br}:=\cap_{x\in Z}\cm^j_x\sset R.
\eeq
 The inclusion $I(Z)^{\bl j\br}\supseteq I(Z)^j$ is obvious, and is often proper. The classical example is  $Z=V(xy,yz,xz)\sset \R^3$.
  Here $xyz\in I(Z)^{\bl 2\br}\smin I(Z)^2$.

The ideal of functions flat on $Z$ is $I(Z)^{\bl\infty\br}=\cap_j I(Z)^{\bl j\br}=\cap_{x\in Z}\cm^\infty_x$.

Denote the completion \wrt the filtration $I(Z)^{\bl\infty\br}$ by $\hR^{(Z)}$, as in equation \eqref{Eq.Whitney.Thm.Completion}.

\item We denote  the partial derivatives by  $\di^{k_1}_{x_1}\dots \di^{k_n}_{x_n}g$ or by $g^{(\uk)}$, using multi-indices.
The condition
\[
\text{``any partial derivative  $\di^{k_1}_{x_1}\dots \di^{k_n}_{x_n}g$, with $\sum_i k_i=|\uk|$ satisfies $|\di^{k_1}_{x_1}\dots\di^{k_n}_{x_n}g|\le ..$"}
\]
 is abbreviated to ``$|g^{(\uk)}|<...$".

The total    $k$'th derivative is a function valued in (dual) $k$-forms, $\cU\stackrel{g^{(k)}}{\to}Hom_\R(Sym^k(\R^n),\R)$.
 Denote its operator norm by $||g^{(k)}||$.
\item
For a closed subset    $Z\sset \cU$   denote by  $C^\infty(\cU,Z)$ the ring of germs of smooth functions at $Z$,
 see example \ref{Ex.Localizations.Germs.along.Z}.
\item
  Fix a filtration by ideals, $R=I_0\supseteq I_1\supseteq \cdots$. Denote $I_\infty:=\cap I_j$.
 Take the corresponding completion, $R\to \hR^{(I_\bullet)}:=\liml_{\leftarrow}\quots{ R }{I_j}$.
 Its elements are (equivalence classes of) Cauchy sequences of functions, $\{f_j\}\in  R $, such that $f_{j+i}-f_j\in I_j$, for all $i,j>0$.
  These elements can be presented also as the formal sums $\sum^\infty_{j=0}g_j\in \prod I_j$.
 These sums are taken up to the equivalence: $\sum g_j\equiv\sum \tg_j$ if $\sum^N(g_j-\tg_j)\in I_N$
 for any $N$. One has the exact sequence:
 \beq
 \bM
0 & \to & \big\{\sum g_j|\ \sum^N_{j=0} g_j\in I_N,\ \forall\ N\big\}
 &\to &\prod^\infty_{j=0} I_j  & \stackrel{\pi}{\to} & \hR^{(I_\bullet)} &\to&0.
\eM \eeq
An element $\sum^\infty_{j=0}g_j\in \hR^{(I_\bullet)}$ is presented by $f\in R$ if  $f-\sum^N_{j=0}g_j\in I_N$  holds  for any $N\ge0$.
\\\parbox{13.2cm}
{  Two filtrations, $I_\bullet$, $\tI_\bullet$ are called equivalent if they satisfy: $I_{k_j}\sseteq \tI_j\sseteq I_{d_j}$,
  where $d_j,k_j\to\infty$ as $j\to \infty$.
 Equivalent filtrations induce isomorphic completions, as on the diagram.}\quad\quad\quad
 $\bM R\to \hR^{(I_\bullet)}\\ \parallel \quad \quad  \parallel\quad \ \\R\to  \hR^{(\tI_\bullet)}\eM$
\eee

\subsection{Whitney extension theorem vs the surjectivity of completion}\label{Sec.Whitney.Extension.Theorem}
 Recall the classical Whitney extension problem:
\beq
\ber\text{Given a closed subset $Z\sset \cU$ and a collection  of
continuous functions, $\{h_\uk\}$, on $Z$,}
\\
\text{ does there exist $f\in C^\infty(\cU)$ with the prescribed partial
derivatives $\{f^{(\uk)}|_Z=h_\uk\}_\uk$?
 }\eer
 \eeq
If $f$ exists then the Taylor expansion of each $h_\uk$ should be expressible via the functions $\{h_{\uk+\tilde\uk}\}_{\tilde\uk}$.

To write this explicitly it is useful to combine,  for each $k$, the partials $\{h_\uk\}_{|\uk|=k}$  into the total derivative:
\beq
Z\stackrel{h_k}{\to}Hom_\R(Sym^k(\R^n),\R).
\eeq
Then the conditions $\{f^{(k)}|_Z=h_k\}_k$ imply the compatibility conditions:
 \beq\label{Eq.Whitney.Compatibility}
\forall\ d,\ k\in \N,\quad  \ \forall\ x,y\in Z:\quad\quad
  h_k|_x-\sum_{0\le \tk\le d} \frac{1}{\tk!} h_{k+\tk}|_y(\underbrace{x-y,\dots,x-y}_\tk)=
  o(||x-y||^d).
 \eeq
The Whitney extension theorem says that these compatibility conditions are also sufficient, \cite{Whitney}.
 Moreover, one can choose $f$ to be real-analytic on $\cU\smin Z$.

\subsubsection{}
A Whitney jet on a closed subset $Z\sset \cU\sseteq\R^n$ is a  collection of functions compatible as in \eqref{Eq.Whitney.Compatibility}, denote it
 by $Z\stackrel{\oplus h_k}{\to}\prod_k Hom_\R(Sym^k(\R^n),\R)$. (Note that we allow infinite sums here.)
Denote by $C^\infty(Z)$ the vector space of all the Whitney jets.
 Define the product on $C^\infty(Z)$ by
 \beq
(\oplus_k h_k)\cdot (\oplus_\tk \th_\tk)=\oplusl_{k,\tk} \bin{k+\tk}{k}\cdot Sym(h_k\otimes h_\tk)\in C^\infty(Z).
\eeq
(Here $Sym(\dots)$ is the total symmetrization of the tensor, this imposes the factor $\bin{k+\tk}{k}$.)

With this product $C^\infty(Z)$ becomes a (commutative, associative, unital) ring. (By the direct check.)
 Whitney extension theorem  is then the exactness of the two  sequences (of rings and ideals):
\beq\label{Eq.Whitney.Thm.Jets}
\bM 0\to I(Z)^{\bl \infty\br}\to& C^\infty(\cU) &\to&  C^\infty(Z)&\to 0
\\& \cup
&& || \\
&  C^\infty(\cU)\cap C^\om(\cU\smin Z)& \to & C^\infty(Z)&\to 0
\eM
\eeq
\bel
Whitney extension theorem identifies $\widehat{C^\infty(\cU)}^{(Z)}\isom{} C^\infty(Z)$, the completion for the filtration $I(Z)^{\bl\bullet\br}$
in equation \eqref{Eq.Whitney.Thm.Completion}. In particular this implies the surjectivity of the completion.
\eel
\bpr
Define the map $\widehat{C^\infty(\cU)}^{(Z)}\stackrel{\Psi}{\to}C^\infty(Z)$ by $\sum g_j\to (\sum g_j|_Z,\sum g^{(1)}_j|_Z,\dots)$.
 If $\sum g_j\equiv \sum \tg_j\in \widehat{C^\infty(\cU)}^{(Z)}$ then $\Psi(\sum g_j)=\Psi(\sum \tg_j)$, see \S\ref{Sec.Notations.Conventions}.v.
 Here each sum $\sum g^{(k)}_j|_Z$ is finite as $g^{(k)}_j|_Z=0$ for $j>k$.
By construction $\Psi$  is a homomorphism of rings. It is injective, as $ker(\Psi)=I(Z)^{\bl \infty\br}\equiv0\in \widehat{C^\infty(\cU)}^{(Z)}$.

For any Whitney jet $\oplus h_k$ take its Whitney representative, $f\in C^\infty(\cU)$. Then its completion, $\hf\in \widehat{C^\infty(\cU)}^{(Z)}$,
 is mapped by $\Psi$ to $\oplus h_k$. Thus $\Psi$ is surjective, hence an isomorphism, and equation \eqref{Eq.Whitney.Thm.Jets} is just
   equation \eqref{Eq.Whitney.Thm.Completion}.
\epr

\

Whitney theorem does not imply the surjectivity for more general
filtrations, even if one assume $V(I_1)=V(I_2)=\cdots$. For example, take a filtration
$C^\infty(\cU)\supset I_\bullet$ such that $\{V(I_j)=Z\}_j$,
$I_j\sset I(Z)^{\bl j\br}$, but $I_\bullet$ is not equivalent to
$I(Z)^{\bl \bullet\br}$. Thus $I_\infty\ssetneq I(Z)^{\bl
\infty\br}$. For any element $\sum g_j\in \hR^{(I_\bullet)}$
Whitney theorem ensures a
 $f\in C^\infty(\cU)$ such that  $f-\sum^N g_j\in I(Z)^{\bl N\br}$, for all $N$. But this $f$ is not a representative of $\sum g_j$ as it does not satisfy
   $f-\sum^N g_j\in I_N$.

Besides, filtrations satisfying  $I_j\sseteq I(Z)^{\bl j\br}$ are rather special, see
  examples \ref{Ex.Various.Filtrations} and \ref{Sec.Examples.surjectivity.of.completions}.
 For the general filtrations the surjectivity question does not seem to be resolvable by Whitney theorem.

\subsubsection{}
We prove that the surjectivity of completion implies Whitney extension theorem for a large class of  sets $Z$.
 Suppose  a closed subset  $Z\sset \cU\sseteq\R^n$ admits a locally finite smooth stratification, i.e. $Z=\coprod^n_{j=0}Z_j$, where each $Z_j$ is the disjoint
 locally finite (possible empty) union  of dimension $j$ manifolds in $\cU$.
 This class of sets is rather large, e.g. any semi-algebraic/sub-analytic set is of this type.

\bel
 For sets with locally finite smooth stratification  the surjectivity of completion  implies the Whitney extension theorem.
\eel
\bpr Take  a Whitney jet  $\{h_k\}_k$ on $Z$. By partition of unity arguments, see e.g. \S\ref{Sec.Preparations.Local.to.Global}, we can pass to
 small open subsets of $\cU$. Thus we assume the finite smooth stratification of $Z$. Moreover, for each connected component of each $Z_j$ we can
  assume that the normal bundle (inside $\cU$) is trivial.
\bee[\bf Step 1.]
\item $Z_0$ is a finite set of points, say $\{p_i\}_i$. For each $p_i$ the restrictions $\{h_k|_{p_i}\}_k$ are the prescribed (total)
 derivatives of $k$'th orders. By the surjectivity of completion (Borel lemma) there exists $f_0\in C^\infty(\cU)$ with these derivatives at $\{p_i\}$.

  Define $\{h_{1,k}:=h_k-f^{(k)}_0|_Z\}_k$, thus $\{h_{1,k}|_{Z_0}=0\}_k$. It remains to resolve the Whitney extension problem for
  $\{h_{1,k}\}_k$ on $Z$.
\item (Induction step.)  Assume we have constructed $f_0,\dots,f_{j-1}$, as before. We have the Whitney jet $\{h_{j,k}\}_k$, these functions
 vanish on $Z_0\cup\dots\cup Z_{j-1}$. Now $Z_j$ is a finite disjoint union of submanifolds of $\cU$, and it is enough to construct $f_j$ on each
  of them separately. (Note that $f_j$ will be flat on the inner boundary of $Z_j$.)

Accordingly, let $Z_j\sset \cU$ be a manifold of $dim=j$. We can assume $Z_j$ is small, in particular we trivialize its normal bundle and split the $\cU$-coordinates
 near $Z_j$: $x_1\dots,x_j$ along $Z_j$, and $y_{j+1},\dots,y_n$ transversal to $Z_j$.
 Take the power series, here $y$ is the multi-variable
 \beq
 \sum^\infty_{k=0} \frac{h_k(\overbrace{y,\dots,y}^k)}{k!}\in C^\infty(\cU_{x_1\dots x_j})[[y_{j+1},\dots,y_n]].
 \eeq
Present it by $f_j\in C^\infty(\cU)$, by the surjectivity of completion. Then $f_j$ extends the Whitney jet $\{h_k\}_k$ on $Z_j$.
 Note that $f_j$ is flat on $\overline{Z_j}\cap(\cup^{j-1}_{i=0} Z_i)$. Define $\{h_{j+1,k}:=h_{j,k}-f^{(k)}_j|_Z\}_k$, thus $\{h_{j+1,k}\}_k$
  all vanish on $\cup^j_{i=0} Z_i$. This finishes the inductive step.
\epr \eee

\beR
By the same argument one can obtain the Whitney extension theorem for the following larger class of sets:
 $Z$ contains a locally finite collection of disjoint manifolds, $\{M_i\}$, whose complement is a discrete, at most countable,
  union of manifolds. Namely, $Z\smin (\coprod M_i)=\coprod N_j$, where $\coprod N_j\sset \coprod \cU(N_j)$, for some disjoint neighbourhoods,
   $\cU(N_j)\cap \cU(N_i)=\empty$ for $i\neq j$.

Sets with locally finite smooth stratifications are obviously of this type.
An  example without a locally finite smooth stratification: $Z$ is a sequence of points with a locally finite set of condensation points.
\eeR

\subsection{Analytic approximation of smooth functions}\label{Sec.Whitney.Approximation}
For a  closed subset   $Z\sset \cU$  take an exhaustion of $\cU\smin Z$ by compactly embedded opens, i.e.
  $\empty\sseteq \cU_1\sseteq \cdots\sset \cU\smin Z$ open bounded such that
 $\overline{\cU_j}\sset \cU_{j+1}$ and $\cup\cU_j=\cU\smin Z$.
\bel\cite[Lemma 6]{Whitney}\label{Thm.Whitney.Approximation}
 For any $\tf\in C^\infty(\cU\smin Z)$ and any sequence of reals $\{\ep_j\to 0\}$
  there exists $f\in C^\om(\cU\smin Z)$ satisfying   $||f^{(k)}-\tf^{(k)}||<\ep_j$  on $\cU\smin (Z\cup \cU_j)$, for   $0\le k\le j<\infty$.
\eel
If $\tf\in C^\infty(\cU)$ then the lemma guarantees: $f\in C^\infty(\cU)\cap C^\om(\cU\smin Z)$.
 Indeed, $\lim_{x\to x_0}f^{(k)}|_x=\tf^{(k)}|_{x_0}$,  for any $k$ and any $x_0\in Z$. Thus each $f^{(k)}$ extends continuously onto the whole $\cU$.

\subsection{Filtered    Whitney theorem of zeros}\label{Sec.Whitney.theorem.of.zeros}
 Whitney theorem on zeros reads: any closed subset $Z\sset \cU$ is presentable as $V(\tau)\sset \cU$ for some $\tau\in C^\infty(\cU)$.
  (In notable difference to the situation Algebraic Geometry/Commutative Algebra.) We give a filtered version, for   $I_\bullet\sset C^\infty(\cU)$.
\bprop\label{Thm.Whitney.Thm.zeros.Strengthened}
 Let $Z\!=\! \overline{\cup_j V(I_j)} \!\sset\! \cU$.
  Then exists $\tau\!\in\! I_\infty\cap C^\om(\cU\smin Z)$ such that $\tau|_Z\!=\!0$
   and $\tau|_{\cU\smin Z}\!>\!0$.
\eprop
\bpr
\bee[\bf 1.]
\item
 For any $j$ we construct $\tau_j\in I_j$ satisfying $V(\tau_j)=V(I_j)$. For any point $p\in \cU\smin V(I_j)$ exists some $\tau_p\in I_j$ such that $\tau_p(p)\neq0$.
 By local compactness this collection $\{\tau_p\}$ can be chosen locally finite on $\cU$. Then define $\tau_j:=\sum(\tau_p)^2\cdot a_p$, where $\{a_p\}$ are
  the relevant non-negative cutoff functions.
\item
 We construct $\ttau\in I_\infty\sset C^\infty(\cU)$ satisfying: $V(\ttau)=Z$ and $\ttau>0$ on $\cU\smin Z$.
  Fix $\{\tau_j\in I_j\}_j$ satisfying $\{V(\tau_j)=V(I_j)\}_j$. We can assume
  $0\le \tau_j\le 1$ and $\tau_j$ are flat on $V(I_j)$. (For example, replace $I_j$ by $I_j\cap I(V(I_j))^{\bl \infty\br}$.)

    Take a decreasing sequence of positive reals $\{\ep_j\to0\}$, and the corresponding small neighbourhoods,
  $\cU_{\ep_j}(Z)$, see \S\ref{Sec.Notations.Conventions}.i.
 Define $\ttau$ on $\{\cU_{\frac{\ep_{j-1}+\ep_j}{2}}(Z)\smin \cU_{\ep_j}(Z) \}_j$ by
     $\ttau(x)=\tau_1(x)\cdots \tau_j(x)$. Extend it in a $C^\infty$ way to $\cU\smin Z$
      such that       $\tau_1 \cdots \tau_{j+1} \le \ttau \le \tau_1 \cdots \tau_j $   on $\cU_{\ep_j}(Z)\smin \cU_{\frac{\ep_j+\ep_{j+1}}{2}}(Z)$.

We have defined a positive function $\ttau\in C^\infty(\cU\smin Z)$.
 It is flat on each $V(I_j)$,  and therefore flat on $Z= \cup_j V(I_j)$.
 Extend $\ttau$ to $Z$ by $0$. Thus $\ttau\in I_j$ for any $j$, and $V(\ttau)=Z$.
\item
  Take a sequence $\{\ep_j\to 0\}$
 and the corresponding bounded opens
 \beq
 \cU_j=\{x|\ |\ttau(x)|>\ep_j\}\cap Ball_{\frac{1}{\ep_j}}(o)\cap \cU.
\eeq
 Thus $\cup_j \cU_j=\cU\smin Z$.
 By   lemma \ref{Thm.Whitney.Approximation} there exists
  $\tau\in C^\infty(\cU)\cap C^\om(\cU\smin Z)$ satisfying $|\tau-\ttau|<\ep_j$ on $\cU_j$. Thus $\tau|_Z=0$, $\tau|_{\cU\smin Z}>0$ and $\tau\in I_\infty$.
\epr
\eee
\bex
\bee[\bf i.]
\item In the trivial case of constant filtration, $\{I_j=I_1\}$, we get: any closed $Z\sset \cU$
  is the zero locus of a smooth function which is real-analytic off $Z$.
\item Suppose the zero sets $\{V(I_j)\}$ stabilize. Then $\cup V(I_j)=V(I_\infty)=V(\tau)$ for some $\tau\in I_\infty$.
 In particular, $I_\infty\neq 0$, regardless of how fast the ideals $I_j$ decrease.
\eee
\eex

\subsection{Persistence of the surjectivity under change of rings}\label{Sec.Persistence.of.surjectivity.under.change.of.rings}

\subsubsection{The local-to-global transition}\label{Sec.Preparations.Local.to.Global}
Suppose we want to represent an element of completion  $\sum g_j\in \widehat{C^\infty(\cU)}^{(I_\bullet)}$  by some $f\in C^\infty(\cU)$.
 \bel
 It is enough to resolve this locally at each point of $\cU$.
\eel
\bpr
Suppose for a cover $\cU=\cup \cU_\al$ the functions $\{f_\al\in C^\infty(\cU_\al)\}_\al$ represent the elements
  $\{\sum g_j|_{\cU_\al}\in \widehat{C^\infty(\cU_\al)}^{(I_\bullet)}\}_\al$.
 Thus
 \beq
 f_\al-\sum^N_{j=0} g_j|_{\cU_\al} \in I_N\cdot C^\infty(\cU_\al),\quad \text{ for each $\al$ and $N$.}
\eeq
 We can assume that the cover is locally finite, by shrinking $\{\cU_\al\}_\al$ if needed.
  Take the corresponding partition of unity,
\beq
\{u_\al\in C^\infty(\cU)\}_\al:\quad\quad
0<u_\al\big|_{ \cU_\al} \le1,\quad\quad\quad\quad
 u_\al\big|_{\cU\smin \cU_\al}=0, \quad\quad\quad\quad
 \sum u_\al=1_\cU.
\eeq
We can assume each $f_\al$ is bounded on $\cU_\al$, thus $u_\al f_\al\in C^\infty(\cU)$.
Define $f:=\sum u_\al f_\al$, then $f\in C^\infty(\cU)$, as the sum is locally finite. And $f$ represents $\sum g_j$:
 \beq
\forall N:\quad \quad  f-\sum^N_{j=0} g_j=\sum_\al u_\al f_\al-\sum^N_{j=0} 1_\cU\cdot g_j=\sum_\al (u_\al f_\al-\sum^N_{j=0} u_\al g_j)\in I_N.
\quad\quad\quad\epr \eeq

\subsubsection{}
 Take a surjective homomorphism of (commutative, associative) rings, $R\stackrel{\phi}{\twoheadrightarrow}\tR$.
 \bel
 \bee[1.]
 \item Suppose the completion of $R$ is surjective, $R\twoheadrightarrow \hR^{(I_\bullet)}$.
  Then the completion of $\tR$ is surjective, $\tR\twoheadrightarrow \hat{\tR}^{(\tI_\bullet)}$, here
    $\tI_\bullet:=\phi(I_\bullet)\cdot \tR$.
\item Suppose the completion of $\tR$ is surjective, $\tR\twoheadrightarrow \hat{\tR}^{(\tI_\bullet)}$.
 Then  the completion of $R$ is surjective, $R\twoheadrightarrow \hR^{(I_\bullet)}$, here $I_\bullet:=\phi^{-1}(\tI_\bullet)$.
 \eee
 \eel
\bpr
\bee[1.]
\item  The homomorphism $R\stackrel{\phi}{\to}\tR$ induces the homomorphism of completions,
  $\liml_\leftarrow\quots{R}{I_\bullet}\stackrel{\hat\phi}{\to}\liml_\leftarrow \quots{\tR}{\tI_\bullet}$.
  Thus we get the commutative diagram with exact rows.
\beq\label{Eq.diagram.inside.proof}
\bM 0\to I_\infty\to R\stackrel{\widehat{}_R}{\to}\hR^{(I_\bullet)}
\\\quad\quad \downarrow \quad \phi\downarrow \quad\quad \downarrow \hat\phi\\
 0\to \tI_\infty\to \tR\stackrel{\widehat{}_\tR}{\to}\widehat\tR^{(\tI_\bullet)} \eM
 \eeq
 If $\phi$ is surjective then $\hat\phi$ too.
  Thus the surjectivity of \ $\widehat{}_R$ implies that of \ $\widehat{}_\tR$.
\item Starting from the filtration $\tI_\bullet\sset \tR$ we get $I_\bullet\sset R$ and the diagram \eqref{Eq.diagram.inside.proof}.
 Then the surjectivity
  of \ $\widehat{}_\tR$ implies that of \ $\widehat{}_R$. \epr
\eee
\bex\label{Ex.Surjectivity.Passage.to.quotient}
\bee[\bf i.]
\item The surjectivity   $C^\infty(\cU)\twoheadrightarrow \widehat{C^\infty(\cU)}^{(I_\bullet)}$ implies 
 $\quots{C^\infty(\cU)}{J}\twoheadrightarrow \widehat{\quots{C^\infty(\cU)}{J}}^{(I_\bullet(mod\ J))}$.
\item The surjectivity    $\quots{C^\infty(\cU)}{J}\twoheadrightarrow \widehat{\quots{C^\infty(\cU)}{J}}^{(I_\bullet(mod\ J))}$
   implies 
    $C^\infty(\cU)\twoheadrightarrow \widehat{C^\infty(\cU)}^{(I_\bullet)}$, provided $I_\infty\supseteq J$.
\eee
\eex

\subsubsection{Localization of a ring} Let $S\sset R$ be a multiplicatively closed set, and $I_\bullet\sset R$ a filtration.
 Take its image under the localization, $I_\bullet[S^{-1}]\sset R[S^{-1}]$.
\bel\label{Thm.Surjectivity.Localization}
Suppose $I_\infty\cap S=\empty$ and for any $q\in S$ exists $\tau_q\in I_\infty$ such that $ 1+\tau_q\in (q)\sset  R$.
  Then $R\twoheadrightarrow \hR^{(I_\bullet)}$
 iff $R[S^{-1}]\twoheadrightarrow \widehat{R[S^{-1}]}^{(I_\bullet[S^{-1}])}$.
\eel
Note that $R\to R[S^{-1}]$ is non-surjective and $I_\bullet\neq R\cap (I_\bullet\cdot [S^{-1}])$. Thus this lemma is not implied by the previous lemma.
\bpr $\Rrightarrow$ Take any element $\sum\frac{g_j}{q_j}\in \widehat{R[S^{-1}]}^{(I_\bullet[S^{-1}])}$, here $g_j\in I_j\sset R$ and $q_j\in S$.
 Choose $\{\tau_j\}$ as in the assumption, then $\sum\frac{g_j}{q_j}$ is equivalent to
  $\sum\frac{g_j(1+\tau_j)}{q_j}\in \hR^{(I_\bullet)}$. Take its representative, $f\in R$. Then the image of $f$ in $R[S^{-1}]$ is the
   representative of $\sum\frac{g_j}{q_j}$.

$\Lleftarrow$ Take an element $\sum g_j\in \hR$. Its image in $\widehat{R[S^{-1}]}^{(I_\bullet[S^{-1}])}$ is presented by some $\frac{f}{q}\in R[S^{-1}]$,
 here $q\in S$. Take $\tau$ as in the assumption, then $\frac{f(1+\tau)}{q}\in R$ presents $\sum g_j$.
\epr
\bex\label{Ex.Localizations.Germs.along.Z}
 Let $R=C^\infty(\cU)$ and $Z=V(I_\infty)$. Suppose no element of the set $S$ vanishes at any point of $Z$. Then
 $C^\infty(\cU)\twoheadrightarrow \widehat{C^\infty(\cU)}^{(I_\bullet)}$ iff
    $C^\infty(\cU)[S^{-1}]\twoheadrightarrow \widehat{C^\infty(\cU)[S^{-1}]}^{(I_\bullet[S^{-1}])}$\!\!.
    Indeed, we have $I_\bullet\!\cap\! S\!=\!\empty$.
     In addition, for any $q\!\in\! S$ take   separating open neighbourhoods of the closed sets,   \S\ref{Sec.Notations.Conventions}.i:
     \beq
     Z\sset \cU(Z),\quad\quad\quad\quad V(q)\sset \cU(V(q)),\quad \quad\quad \cU(Z)\cap  \cU(V(q))=\empty.
     \eeq
Then take $\tau\in I_\infty$      such that $\tau|_{\cU(V(q))}=-1$.
 It exists e.g. by proposition \ref{Thm.Whitney.Thm.zeros.Strengthened}. Then $1+\tau\in (q)$, now invoke lemma \ref{Thm.Surjectivity.Localization}.

Some particular cases:
\bei
\item For a point $Z=\{x\}\sset  \cU$, and $S=\{$all the functions that do not vanish at $x\}$, we get the germs at the point,
  $C^\infty(\cU)[S^{-1}]=C^\infty(\cU,x)$.
\item For a closed subset $Z\sset\cU$ and  $S=\{$all the functions that do not vanish at any point of $Z\}$, we get the germs along $Z$,
  $C^\infty(\cU)[S^{-1}]=C^\infty(\cU,Z)$.
\item Let $Z\times\cV\sset \cU\times\cV$, for a closed $Z\sset \cU$. Take $S=\{$all the functions that do not vanish at any point of $Z\times\cV\}$.
 Then   $C^\infty(\cU\times\cV)[S^{-1}]=C^\infty\big((\cU,Z)\times\cV\big)$.
\eei
\eex

Thus the surjectivity question for rings of germs, e.g. $C^\infty(\cU,Z)$, is reduced to the ring $C^\infty(\cU)$.

\beR\label{Re.surjectivity.of.completions.simple.reduction}
 For some rings/filtrations the surjectivity of completion   can be proved by the following simple
argument. Assume $(R,\cm)$ is local and $I_j\sseteq \cm^{d_j}$, with
$d_j\to \infty$.   Then the completion $R\to \hR^{(\cm)}$ factorizes
through
 $R\to \hR^{(I_\bullet)}\to \hR^{(\cm)}$. Thus, if the map $R\to \hR^{(\cm)}$ is surjective and the map  $\hR^{(I_\bullet)}\to\hR^{(\cm)}$
  is injective, the map $R\to \hR^{(I_\bullet)}$ is surjective.
 However, a necessary condition for the injectivity   $\hR^{(I_\bullet)}\to\hR^{(\cm)}$ is $\cap I_j\supseteq\cm^\infty$.
  And this does not hold for many filtrations. Thus one must prove the surjectivity separately.
\eeR

\subsection{Surjectivity for rings of functions on manifolds and other spaces}\label{Sec.Surj.For.manifolds.and.schemes}
 Let $\cU\sseteq\R^n$ and $X\sset \cU$ a closed subset, with the defining ideal $I(X)$.
  Take the ring of  Whitney jets, $C^\infty(X)$, see \S\ref{Sec.Whitney.Extension.Theorem}. By Whitney theorem we
   identify $C^\infty(X)=\quots{C^\infty(\cU)}{I(X)^{\bl \infty\br}}$.
 For sets with dense interior, $X=\overline{Int(X)}$, the elements of $C^\infty(X)$ are those elements of
 $C^\infty(Int(X))$ whose derivatives extend continuously onto $X$.
 In this case $I(X)^{\bl \infty\br}=I(X)$ and $C^\infty(X)=\quots{C^\infty(\cU)}{I(X)}$.

    If $X\sset \cU$ is a (closed) $C^\infty$-submanifold then the ring of smooth functions on $X$ can be presented as $\quots{C^\infty(\cU)}{I(X)}$.
 For a $C^\infty$-algebraic scheme this is the ring of $C^\infty$-regular functions, \cite{Joyce}.
 For  analytic subset  $X\sset \cU$    the ring of real-analytic functions is $C^\om(X):=\quots{C^\om(\cU)}{I(X)}$.

In theorems \ref{Thm.Completion.Nec.Condition} and \ref{Thm.Completion.Suff.Condition} we address the surjectivity for the ring $\quots{C^\infty(\cU)}{J}$.
 In particular this gives the necessary and sufficient condition for the surjectivity  $C^\infty(X)\twoheadrightarrow \widehat{C^\infty(X)}^{(I_\bullet)}$
  in all these cases.

\section{The necessary condition for the surjectivity of completion map}\label{Sec.Surjectivity.Nec.Condition}
Let $R=\quots{C^\infty(\cU)}{J}$, for some open $\cU\sseteq\R^n$  and some ideal $(0)\sseteq J\sset C^\infty(\cU)$. Fix a filtration $I_\bullet\sset R$,
 and denote $Z:=V(I_\infty)$.
For an open subset $\cU_0\sset \cU$ we take the restriction  $I_\bullet\big|_{\cU_0}\sset R|_{\cU_0}:=\quots{C^\infty(\cU_0)}{C^\infty(\cU_0)\cdot J}$.

 In this section we prove:
\bthe\label{Thm.Completion.Nec.Condition}
  Suppose the completion map is surjective, $R \twoheadrightarrow \hR^{(I_\bullet)}$.
Take any compactly embedded open subset,  $\cU_0\sset \overline{\cU_0}\sset \cU$, where $\overline{\cU_0}$ is compact in
$\R^n$. The filtration restricted to $\cU_0$ is equivalent to a particular form:
 $I_\bullet|_{\cU_0}\sim I_\infty+I_\bullet\cap I(Z)^{\bl\bullet\br}|_{\cU_0}$.
\ethe

\

First we introduce the   loci $Z_i\sset \cU$ where the ideals $I_\bullet$ have prescribed orders \wrt the maximal ideals.
  The surjectivity of $R\to \widehat{R}\ \!^{(I_\bullet)}$
 imposes heavy restrictions on $I_\bullet$ along these loci, lemma \ref{Thm.Zero.Loci.Stabilization}.
  Using these restrictions we prove theorem \ref{Thm.Completion.Nec.Condition}.

As a corollary we get the necessary condition for localized rings, \S\ref{Sec.Surjectivity.Nec.Criterion.Localization}.

\subsection{The loci of prescribed orders}
 Take the zero locus $V(J)\sseteq \cU$. Define the loci of $i$'th order:
 $V_i(I_j):=\{x|\ I_j\sseteq \cm^i_x\sset R\}\sset V(J)\sseteq\cU$. Define their limiting loci:
 \beq\label{Def.multiple.chain.of.loci}
 V_i(I_1)\sseteq V_i(I_2)\sseteq \cdots\sseteq  \cupl_j V_i(I_j)=:Z_i.
 \eeq
These loci $\{Z_i\}$ satisfy $V(J)\supseteq Z_1\supseteq Z_2\supseteq \cdots$, they are not necessarily closed.
 Define $Z_\infty=\cap Z_i$.

\bex\label{Ex.Various.Filtrations} \bee[\bf i.]
\item Take the filtration $\{I_j:=I^j\}$. Then $Z_1=\cdots=Z_\infty =V(I)\sseteq V(J)$.
\item Let $\{X_j\}$ be an infinite collection of   compact manifolds in $\R^n$.
 Take their defining ideals, $\{I(X_j)\}_j$.
 Take the filtration
 \beq
 I_j:=I(X_1)\cap I(X_2)^2\cap\cdots\cap I(X_j)^j.
\eeq
  Then \
 $Z_1=\cupl_{j\ge1}X_j\supset  Z_2=\cupl_{j\ge2}X_j\supset\cdots\supset  Z_i=\cupl_{j\ge i}X_j\supset\cdots$.
\eee \eex

With no assumptions on $I_\bullet$ (or on the surjectivity $R\twoheadrightarrow\hR$), the chains in \eqref{Def.multiple.chain.of.loci} do not necessarily stabilize, and the loci $\{Z_i\}$
 are not necessarily  closed, see the last example.

\subsection{Stabilization of the loci of prescribed orders}

\bel\label{Thm.Zero.Loci.Stabilization}
Let $R=\quots{C^\infty(\cU)}{J}$ and suppose the completion map is surjective, $R \twoheadrightarrow \hR^{(I_\bullet)}$.
\bee[1.]
\item  For any compactly embedded open subset $\cU_0\sset\overline{\cU_0}\sset \cU$
   the restricted filtration $I_\bullet\big|_{\cU_0}$ is equivalent to the filtration
   \[\{I_j\cap I(Z_1)\cap I(Z_2)^{\bl 2\br}\cap\dots\cap I(Z_j)^{\bl j\br}\big|_{\cU_0}\}.
   \]
\item     In particular,
  the restrictions of chains in  \eqref{Def.multiple.chain.of.loci} stabilize,   $Z_i\cap \cU_0 =V_i(I_j)\cap \cU_0$ for $j\gg i$.
   Therefore   $\{Z_i\sset \cU\}$ are closed for all $1\le i\le \infty$.
\eee
\eel
\bpr
 To check the equivalence of the filtrations it is enough to show: for any $j$ exists $d_j<\infty$ such that
 \beq
 I_{d_j}\big|_{\cU_0}\sseteq I_j\cap I(Z_1)\cap I(Z_2)^{\bl 2\br}\cap\dots\cap I(Z_j)^{\bl j\br}\big|_{\cU_0}\sset
   R|_{\cU_0}=\quots{C^\infty(\cU_0)}{J\cdot C^\infty(\cU_0)}.
 \eeq
  For this it is enough to show: for any $i$ the chain $\{V_i(I_j)\big|_{\cU_0}\}_j$ in   \eqref{Def.multiple.chain.of.loci} stabilizes.
   We prove this by induction on $i$.

{\bf Case $i=1$, for $Z_1$.} Suppose the loci $V_1(I_j)\big|_{\cU_0}\sset  V(J)\cap \cU_0$ do
not stabilize, i.e. $I_j\not\sseteq I(Z_1)$, for any $j$. Replace $I_\bullet$ by its (equivalent) subsequence
that satisfies:
 $V_1(I_1)\big|_{\cU_0}\ssetneq V_1(I_2)\big|_{\cU_0}\ssetneq \cdots$. Fix a sequence of points,
 $p_{j}\in V_1(I_{j+1})\smin V_1(I_j)\big|_{\cU_0}\sset V(J)\cap \cU_0$,
   and small balls, $Ball_{\ep_{j}}(p_{j})\cap V_1(I_j)=\empty$. Accordingly fix a sequence of elements, $\{g_j\in I_j\}$, and the corresponding functions
  $\{\tg_j\in C^\infty(\cU)\}$,  satisfying:
\beq \tg_j(x)=\Bigg\{\ber 0,\ x\not\in Ball_{\ep_{j}}(p_{j})\\>0,\
x\in Ball_{\ep_{j}}(p_{j})\\j!,\ x=p_j
\eer \eeq
(These are constructed from elements of $I_j$ by using the standard bump functions.)

Suppose the element
$\sum^\infty_{j=0} g_j\in \hR^{(I_\bullet)}$ is presented by some
$f\in R$, see \S\ref{Sec.Notations.Conventions}.vi.
 Then $f$ must be presented by $\tf\in C^\infty(\cU)$, which is
 bounded on the compact set $ \overline{\cU_0}$. But this contradicts the construction: $\tf(p_{N-1})=\sum^N_{j=0}\tg_j(p_{N-1})\ge N!$.

{\bf The case of arbitrary $i$} is similar. Assuming the statement for $i$,
i.e. for $Z_1,\dots,Z_i$, we pass to the (equivalent)  filtration
satisfying:
\beq Z_1\cap \cU_0 =V_1(I_1) \cap \cU_0 ,\quad\quad
Z_2\cap \cU_0 =V_2(I_1) \cap  \cU_0,\quad\quad \dots\quad\quad
Z_i\cap \cU_0=V_i(I_1)\cap \cU_0.
\eeq
Suppose the loci
$\{V_{i+1}(I_j)\cap \cU_0\}_j$ do not stabilize, i.e. $I_j\not\sseteq I(Z_{i+1})^{\bl i\br}$ for any $j$.
As before, we can assume
$V_{i+1}(I_1)\cap \cU_0 \ssetneq V_{i+1}(I_2)\cap \cU_0 \ssetneq \cdots$.
 As before, fix a sequence of points and  small balls:
\beq
 p_j\in V_{i+1}(I_{j+1})\smin V_{i+1}(I_j),\quad\quad\quad\quad
 Ball_{\ep_j}(p_j)\cap V_{i+1}(I_j)=\empty.
 \eeq
  Fix a sequence of elements $\{g_j\in I_j\}$ satisfying $ord_{p_j}g_j=i$. Take their representatives $\{\tg_j\in C^\infty(\cU)\}$ satisfying
    $ord_{p_j}\tg_j=i$.
  By using bump functions we can assume $\tg_j|_{\cU\smin Ball_{\ep_j}(p_j)}=0$.
   Moreover,   we can assume the Taylor expansion $\tg_j(x)=Q_{i}(x-p_j)+O(x-p_j)^{i+1}$, where $Q_i$ is a homogeneous polynomial of degree $i$,
  with large coefficients.
 More precisely, as $J\not\supseteq\cm^{i}_{p_j}$, we have the vector subspace    of positive dimension:
 \beq
 \quot{\cm^i_{p_j}+J}{\cm^{i+1}_{p_j}+J}\sset \quots{R}{\cm^{i+1}_{p_j}}.
 \eeq
 Thus we choose $0\neq [Q_i]\in \quots{\cm^i_{p_j}+J}{\cm^{i+1}_{p_j}+J}$ with large coefficients, and then
 $\tg_j(x)+q$ has large $i$'th derivative for any $q\in J$.

Finally, if $\sum g_j\in \hR$ is presented by some $f\in R$ then  any representative  $\tf\in C^\infty(\cU)$ satisfies:
\beq
\sum^N_{j=1} \tg_j-\tf\in \tI_N+J,\quad \text{ for any } N.
\eeq
 Here the derivative $\tf^{(i)}$ must be bounded on $\cU_0$. But, by construction, its values on the sequence of points $\{p_i\}$
 are unbounded.
\epr
\bex\label{Ex.Ideal.Polynomially.Generated}
This lemma forbids numerous filtrations, even with polynomially generated ideals.
\bee[\bf i.]
\item
Fix a sequence $\{p_i\}$ of distinct points
 in $\R^1$, converging to $0$.
 Let
 \beq
  I_j=\bl\prod^j_{i=1}(x-p_i)^i\br\sset C^\infty(\R^1)=:R.
  \eeq
 Here none of $V_i(I_j)$ stabilizes in $j$.
   Thus the completion map $R\to \hR^{(I_\bullet)}$ is not surjective.
\item   If $I_\infty=(0)\sset R$ then the completion map is not surjective. Indeed, the surjectivity forces $V_1(I_j)$
to stabilize (lemma \ref{Thm.Zero.Loci.Stabilization}), and then $V_1(I_\infty)=V_1(I_j)\ssetneq V(J)$,
 by \S\ref{Sec.Whitney.theorem.of.zeros}. This goes in notable difference to the Noetherian rings Commutative Algebra, where the non-stabilizing
  filtrations usually satisfy  $I_\infty=(0)$.
\eee
\eex

\beR The restriction to the compactly embedded subsets,
$(\dots)\big|_{\cU_0}$,  is important. For example, let $\cU=(0,1)\sset
\R^1$
 and $I_j=\big\{f|\ f=0\ on\ [\frac{1}{j},1-\frac{1}{j}]\big\}\sset C^\infty(0,1)$.
   The chain $V_i(I_j)\sset V_i(I_{j+1})\sset\cdots$ does not stabilize, for any $i$.
 But the completion is surjective, $C^\infty(0,1)\twoheadrightarrow \widehat{C^\infty(0,1)}\ \!^{(I_\bullet)}$.
 Indeed, for any sequence $\{g_j\in I_j\}$ we have   $\sum g_j\in C^\infty(0,1)$, as the sum   is locally finite on $(0,1)$.
\eeR

\beR
The sequence of subsets $Z_1\supseteq Z_2\supseteq\cdots$ of \eqref{Def.multiple.chain.of.loci} does not necessarily stabilize.
For example, take a converging sequence $\R^n\ni p_i\to p_0$.
 Take
 \beq
 Z_i:=\{p_0,p_i,p_{i+1},p_{i+2},\dots\}\sset \R^n,\quad \quad\quad
  I_j:=I(Z_1)\cap \cdots\cap I(Z_j)^{\bl j\br}\sset R:=C^\infty(\R^n).
  \eeq
   The sequence $\{Z_j\}$ does not stabilize, but the completion is surjective,
  $R\twoheadrightarrow \hR^{(I_\bullet)}$. Indeed, fix $\sum g_j\in \hR^{(I_\bullet)}$, here $g_j\in I_j$. This element is presented by $f\in R$ iff $f$
   has the fixed jets, $jet_{i,p_i}f=jet_{i,p_i}\sum^i_{j=1}g_j$. By Whitney extension,  \S\ref{Sec.Whitney.Extension.Theorem},
    we recover $f$ with such derivatives.
   The only Whitney compatibility condition to be checked is at $p_0$, this reads: each $g_k$ is flat at $0$. And this holds as each $g_k$
    vanishes on the whole sequence $\{p_i\}$.
\eeR

\subsection{Proof of theorem \ref{Thm.Completion.Nec.Condition}}
Take the completion, $R\stackrel{\pi}{\to}\hR^{(I_\bullet)}$.
\bee[\bf Step 1.]
\item
We claim: $\cap_j(\pi (I_j)\!\cdot\! \hR^{(I_\bullet)})\!=\!(0)$. 
Indeed, for any $N$ the ideal $\pi (I_N)\!\cdot\! \hR^{(I_\bullet)}$ is sent to
  $(0)\sset \quots{\hR^{(I_\bullet)}}{\pi (I_N)\cdot \hR^{(I_\bullet)}}   \cong\quots{R}{I_N}$.
 Thus we get the projective limit of the zero ideals,
 \beq
 \cap_j(\pi (I_j)\cdot \hR^{(I_\bullet)})=\liml_{\leftarrow}(0)=(0)\sset \liml_{\leftarrow}\quots{R}{I_\bullet}.
\eeq

\item
Consider the quotient ring $\quots{R}{I_\infty}$, with its filtration $\quots{I_\bullet}{I_\infty}$.
For any fixed $i\ge1$ and any point $p\in V(I_\infty)\sseteq\cU$ the quotient filtration
\beq
\tI_\bullet:=\quot{I_\bullet+\cm^i_p}{I_\infty+\cm^i_p}\sset \quot{R}{I_\infty+\cm^i_p}
\eeq
stabilizes. (As a descending chain of ideals in the Artinian ring.)
 Thus we can identify $\quots{R}{I_\infty+\cm^i_p}$ with its $\tI_\bullet$-completion, and identify the  $\tI_\bullet$-filtration with its completed version:
\beq\bM
\quot{R}{I_\infty+\cm^i_p}&\isom{}&\widehat{\quots{R}{I_\infty+\cm^i_p}}^{(\tI_\bullet)}&\isom{}&
\quot{\hR^{(I_\bullet)}}{\pi(I_\infty)  \hR^{(I_\bullet)}\!+\!\pi(\cm^i_p)  \hR^{(I_\bullet)}}
\\
\cup&&\cup&&\cup
\\
\tI_\bullet&\isom{}&\widehat{\quots{R}{I_\infty+\cm^i_p}}^{(\tI_\bullet)}\cdot \tI_\bullet&\isom{}&
\quot{\pi(I_\bullet)  \hR^{(I_\bullet)}\!+\!\pi(\cm^i_p)  \hR^{(I_\bullet)}}{\pi(I_\infty)  \hR^{(I_\bullet)}\!+\!\pi(\cm^i_p)  \hR^{(I_\bullet)}}
\eM.\eeq


By step 1 the bottom-right corner of this diagram converges to zero, for $\bullet\to \infty$. Therefore  $\tI_\bullet$ stabilizes to zero,
 for any $p\in V(I_\infty)$ and any $i$. Therefore $V(I_\infty)=\cap Z_i(\quots{I_\bullet}{I_\infty})=Z_\infty(\quots{I_\bullet}{I_\infty})$.

But the surjectivity of the initial completion, $R  \twoheadrightarrow   \hR^{(I_\bullet)}$, implies the surjectivity of the completion of the quotient,
 $\quots{R}{I_\infty}  \twoheadrightarrow   \quots{\hR^{(I_\bullet)}}{\hI_\infty}$, see example \ref{Ex.Surjectivity.Passage.to.quotient}.
  Thus, by lemma \ref{Thm.Zero.Loci.Stabilization}, we have equivalence of the filtrations,
   $\quots{I_\bullet}{I_\infty}|_{\cU_0}\sim \quots{I_\bullet}{I_\infty}\cap\quots{I(Z)^{\bl \bullet\br}}{I_\infty}\!|_{\cU_0}$.
   Thus $I_\bullet\!|_{\cU_0}\sim I_\infty+I_\bullet \cap I(Z)^{\bl \bullet\br} |_{\cU_0}$, for
    $Z=Z_\infty(\quots{I_\bullet}{I_\infty})=V(I_\infty)$. \epr
   \eee

\beR
One would like to strengthen theorem \ref{Thm.Completion.Nec.Condition} and to obtain further special properties of $I_\bullet$.
 For example, suppose all $I_\bullet$
 are analytically generated, how does the surjectivity restrict the primary decomposition of $I_j$?
 As $V(I_\bullet)$ stabilizes, the minimal primes stabilize, and the corresponding primary ideals form a filtration. But not much can be said
  about the embedded primes, because of the sufficient condition, theorem \ref{Thm.Completion.Suff.Condition}.
   For example, let $\{\cq_{j,k}\}_k$ be the primary ideals corresponding to the minimal primes of $I_j$. Then $\{\cq_{j,k}\}_j$ is a filtration of
    $\cq_{1,k}$, for each $k$. Assume $\cq_{j,k}\sseteq \cq_{1,k}^j$. Then, for any sequence $\{emb_j\}$, corresponding to embedded primes, the filtration
     $ emb_\bullet\cap (\cap_k \cq_{\bullet,k})$ induces the surjective completion map.
\eeR

\subsection{The necessary criterion for localizations}\label{Sec.Surjectivity.Nec.Criterion.Localization}
Let $R=\quots{C^\infty(\cU)}{J}$  and  $S\sset R$ be a multiplicatively closed set. Assume no element of $S$ vanishes at any point of $Z$. Take a filtration
 $I_\bullet\sset R[S^{-1}]$. For any open subset define the restriction $I_\bullet|_{\cU_0}\sset R[S^{-1}]|_{\cU_0}=\quots{C^\infty(\cU_0)}{C^\infty(\cU_0)J}[S^{-1}]$.
\bcor
If $R[S^{-1}]\twoheadrightarrow \widehat{R[S^{-1}]}^{(I_\bullet)}$ then
 $I_\bullet|_{\cU_0}\sim I_\infty|_{\cU_0}  +I_\bullet\cap I(Z)^{\bl \bullet\br}|_{\cU_0}$ for any compactly embedded
  open subset $\cU_0\sset \overline{\cU_0}\sset \cU$.
\ecor
\bpr
For $R\stackrel{\phi}{\to} R[S^{-1}]$ take $\phi^{-1}(I_\bullet)\sset R$. We get
 $R\twoheadrightarrow \hR^{(\phi^{-1}(I_\bullet))}$, by example \ref{Ex.Localizations.Germs.along.Z}. Thus (theorem \ref{Thm.Completion.Nec.Condition})
 \beq
 \phi^{-1}(I_\bullet)|_{\cU_0}\sim \cap_j \phi^{-1}(I_j)|_{\cU_0}+ \phi^{-1}(I_\bullet)\cap  I(Z)^{\bl \bullet\br} |_{\cU_0}.
 \eeq
 For any ideal $\ca\sset R[S^{-1}]$ one has $\ca=R[S^{-1}]\cdot\phi\big(\phi^{-1}(\ca)\big)$. Indeed, take any system of generators $\{a_i\}$ of $\ca$,
  not necessarily finite. Clear the denominators, thus can assume $\{a_j\in \phi(\ca)\}$. Thus $\phi(\phi^{-1}(\ca))=\phi(R)\{a_j\}$.
   Therefore  $I_j=R \phi(\phi^{-1}(I_j))|_{\cU_0}\sseteq I_\infty+ I_{d_j}\cap  I(Z)^{\bl d_j\br} |_{\cU_0}$.
\epr

\section{The sufficient condition for the surjectivity}\label{Sec.Surjectivity.proof.of.Suffic.Cond}
By \S\ref{Sec.Persistence.of.surjectivity.under.change.of.rings} the question is reduced to  the particular ring $C^\infty(\cU)$.
 Take a filtration $I_\bullet\sset R:=C^\infty(\cU)$, let $Z=V(I_\infty)\sset \cU$.
 For each point   $p\in \cU$ we take the ring of germs, $C^\infty(\cU,p)$, see example \ref{Ex.Localizations.Germs.along.Z}, and its (localized) filtration.
\bthe\label{Thm.Completion.Suff.Condition}
Suppose for each point  $p\in \cU$ the localized filtration $I_\bullet|_{(\cU,p)}$
  is (locally) equivalent to a filtration of type $ \ca+\sum_k  \cb_{k,\bullet}$, where
\bei
\item The ideal  $\ca$  does not depend on $\bullet$;  the  summation over $k$ is finite.
\item The ideals $\{\cb_{k,\bullet}\}$ satisfy: $\cc_k\cdot \cb_{k,\infty}\sseteq \cb_{k,j}\sseteq \cc_k\cap \cb_k^{\bl d_j\br}$, where the ideals
 $\cc_k$, $\cb_k$ do not depend on $j$; the zero loci are constant, $V(\cb_{k,\infty})=V(\cb_{k,j})$; and $d_j\to\infty$.
\eei

Then the two sequences are exact:
\beq\label{Eq.Completion.Surjective.Theorem}
\bM 0\to I_\infty \to& R &\to& \hR^{(I_\bullet)}&\to 0
\\& \cup
&& || \ \ &
\\
&  R\cap C^\om(\cU\smin Z)& \to &\hR^{(I_\bullet)}&\to 0
\eM
\eeq

Moreover, for any neighbourhood $Z\sset \cU(Z)$ and any element  $\hf\in \hR^{(I_\bullet)}$ there exists a representative
 $f\in  C^\infty(\cU)\cap C^\om(\cU\smin Z)$ that satisfies  $f>0$ on $\cU\smin \cU(Z)$.
\ethe
 \bpr
  Given $\sum g_j\in \hR^{(I_\bullet)}$, with $g_j\in I_j\sset C^\infty(\cU)$, we should construct its representation  $f\in
C^\infty(\cU)$,  satisfying $f-\suml^N_{j=0}g_j\in I_N$ for any $N$.
 There are several steps:
\bee[\em Step 1.]
\item We  reduce the proof to the ring $C^\infty(Ball_1(o))$ and a very particular filtration, equation \eqref{Eq.particular.filtration}.
\item   We bound the growth of derivatives of $g_j$ on the locus $V(\cb)\sset Ball_1(o)$.
\item We construct $f$ from $\{g_j\}$ using the cutoff functions with controlled growth.
\item We verify that $f$ represents $\sum g_j$ in $\hR^{(I_\bullet)}$. This proves exactness of the first row of  \eqref{Eq.Completion.Surjective.Theorem}.
\item For this $f$ we use Whitney approximation theorem to establish the second row of  \eqref{Eq.Completion.Surjective.Theorem}.
\eee

\bee[\bf Step 1.]
\item  (Simplifying the filtration)
By the partition of unity argument, \S\ref{Sec.Preparations.Local.to.Global}, we can pass to the ring $R:=C^\infty(Ball_1(o))$
 and replace $I_\bullet$ by the equivalent filtration as in the
assumptions:
\beq
I_\bullet=\ca+\sum_k  \cb_{k,\bullet}.
\eeq

  An element of $\hR$ is $a_0+\sum_{j\ge1} (a_j+\sum_k b_{k,j})$, where $a_j\in \ca$ and
   $b_{k,j} \in   \cb_{k,j}$.
   We should construct a representative $f\in R$ that satisfies: $f-a_0-\sum^N_{j=1}(a_j+\sum_k b_{k,j})\in I_N$ for any $N$.
   As $\ca\sseteq I_j$  for each $j$, one can omit $\{a_j\}$.

    Moreover, as the summation $\sum_k$ is finite, it is enough to find a representative for each $\sum_{j\ge1}  b_{k,j}$, with $k$ fixed.
 This reduces the statement to the filtration $\{ \cb_j\}_j$, with $\cc\cdot \cb_\infty\sseteq \cb_j\sseteq \cc\cap \cb^{\bl d_j\br}$,
  $V(\cb_\infty)=V(\cb)$.
Finally, pass to an equivalent filtration with $d_j=j$.

Summarizing, it is enough to establish the surjectivity
$R\twoheadrightarrow \hR^{(I_\bullet)}$ for $R=C^\infty(Ball_1(o))$ and the  filtration of the
particular type:
\beq\label{Eq.particular.filtration}
 \cc\cdot \cb_\infty\sseteq I_j\sseteq \cc\cap \cb^{\bl j\br},\quad \text{ with } V(\cb_\infty)=V(\cb_j).
\eeq

\item
We have $\{g_j\in I_j\}$ for the specific filtration of the ring $C^\infty(Ball_1(o))$ as in \eqref{Eq.particular.filtration}.
 By slightly shrinking the ball we can assume $g_j\in C^\infty(\overline{Ball_1(o)})$, in particular each derivative of each $g_j$ is bounded.
 Take the zero locus $Z_\cb:=V(\cb)\sset Ball_1(o)$.

 For any $0\le k<j<\infty$ and any $x\in Ball_1(o)$ we bound the operator norm of $k$'th derivative:
\beq\label{Eq.bound.on.derivatives.g}
  ||g^{(k)}_j|_x||\le C_{g_j}\cdot dist(x,Z_\cb)^{j-k}.
\eeq
(Here $\{C_{g_j}\}$ are some constants that depend on $g_j$.)

Indeed, for any $x\in Ball_1(o)\smin Z_\cb$ take some $z\in Z_\cb$ for which $dist(x,z)-dist(x,Z_\cb)=\ep\ll dist(x,Z_\cb)$.
 By the assumption $g_j\in \cm_z^j$, thus $g^{(k)}_j|_z=0$ for $k<j$. Therefore the Taylor expansion  (in $Ball_{dist(x,z)}(z)$) with remainder
  gives:
\beq
g_j( x )=\frac{j}{j!}\cdot \intl^1_0(1-t)^{j-1} g^{(j)}_j\big|_{(\uz+t(\ux-\uz))}(\underbrace{x-z,\dots,x-z}_{j})dt.
\eeq
(See \S\ref{Sec.Notations.Conventions}.iii.)
     Note the bounds
\beq
|g^{(j)}_j\big|_{(\uz+t(\ux-\uz))}(\underbrace{x\!-\!z,\dots,x\!-\!z}_{j})|\!\le \!||g^{(j)}_j\big|_{(\uz+t(\ux-\uz))}||\cdot ||x-z||^j,
\quad \quad
||x -z||^j\le \big(dist( x ,Z_\cb)+\ep\big)^j.
\eeq
The derivatives   $||g^{(j)}_j||$ are bounded on  $Ball_1(o)$. Thus
  $|g_j( x )|\le C_0\cdot dist( x ,Z_\cb)^j$, for a constant $C_0$.

The bounds on the derivatives,  $||g^{(k)}_j|_x||\le\dots$, are obtained
in the same way, by Taylor expanding $g^{(k)}_j$ at $z$.
 For each $k$ we get its constant $C_k$. Then define $C_{g_j}:=\maxl_{k<j}(C_k)$.

\item
 We use a particular cutoff function with controlled growth of total derivatives:

{\bf Theorem 1.4.2 of \cite[pg. 25]{Hormander}} {\em  For any
compact set with its neighbourhood, $Z\sset \cU\sset\R^n$,  and
 a positive decreasing sequence  $\{d_j\}$, satisfying $\sum d_j< dist(Z,\di\cU)$, there exists a smaller neighbourhood, $Z\sset \cU(Z)\ssetneq \cU$,  and
     a    function $\tau\in C^\infty(\R^n)$ satisfying
\bee
\item $\tau|_{\R^n\smin \cU}=0$, $\tau|_{\cU(Z)}=1$;
\item for any $k$  and for any $x\in \cU$  the norm of $k$'th derivative is bounded: $||\tau^{(k)}|_x||\le \frac{C^{k}\cdot }{d_1\cdots d_k}$.
\eee
}
(Here   the constant $C$ depends only on the dimension $n$.)

\

In our case the subset  $Z_\cb\sset  Ball_1(o) $ is closed and we can assume it is compact by shrinking the ball.
  Take the  neighbourhood,   $\cU_\ep(Z_\cb)$, see \S\ref{Sec.Notations.Conventions}.i.
 Fix a  decreasing sequence of positive numbers $\{\ep_j\}$,   $\ep_j\to 0$.
  Assume it decreases fast, so that for each $j$ exists a cutoff function satisfying:
 \beq\label{Eq.cutoff.func.derivatives.bounded.growth}
 \tau_j|_{\cU_{\ep_{j+1}}(Z)}=1,\quad\quad
   \tau_j|_{Ball_1(o)\smin \cU_{\ep_j}(Z)}=0,\quad\quad
   \text{and $||\tau^{(k)}_j||$ is bounded as before, for any $k$.}
 \eeq
Define $f(x):=\sum_j  \tau_j(x)\cdot g_j(x)$. We claim that $f\in C^\infty(Ball_1(o))$, provided the sequence $\{\ep_j\}$ decreases quickly.

 The statement $f\in C^\infty(\R^n\smin Z_\cb)$ is obvious, as for any $x\in \R^n\smin Z_\cb$ the summation is finite.
 To check the behaviour on/near $Z_\cb$ we bound $\uk$'th partial derivatives:
\begin{multline}
\Big|\big(\tau_j(x)\cdot g_j(x)\big)^{(\uk)}\Big|=
\Big|\suml_{0\le \ul\le \uk}\bin{|\uk|}{\ul}\tau^{(\ul)}_j(x)\cdot g^{(\uk-\ul)}_j(x)\Big|\stackrel{eq.\ \eqref{Eq.bound.on.derivatives.g}}{\le}
\\
\suml_{0\le \ul\le \uk}\bin{|\uk|}{\ul}\Big|\tau^{(\ul)}_j(x)\cdot C_{g_j}\cdot dist(x,Z_\cb)^{j-|\uk|+|\ul|}\Big|
\stackrel{eq.\ \eqref{Eq.cutoff.func.derivatives.bounded.growth}}{<}
\\C_{g_j}\cdot\suml_{0\le \ul\le
\uk}\bin{|\uk|}{\ul}C^{|\ul|} \cdot
\frac{dist(x,Z_\cb)^{j-|\uk|+|\ul|}}{d_1\cdots d_{|\ul|}}
<
 C_{g_j}\cdot dist(x,Z_\cb)\cdot \suml_{0\le \ul \le \uk}\bin{|\uk|}{\ul}\frac{C^{|\ul|}  }{d_1\cdots
d_{|\ul|}} \ep_j^{j-|\uk|+|\ul|-1}.
\end{multline}
We assume the sequence $\{\ep_j\}$ decreases fast to ensure, for $j>|\uk|+1$:
\beq
C_{g_j} \cdot \suml_{0\le \ul\le \uk}\bin{|\uk|}{\ul}\frac{C^{|\ul|}  }{d_1\cdots d_{|\ul|}}\ep_j^{j-|\uk|+|\ul|-1}<\frac{1}{j!}.
\eeq
Present $f^{(\uk)}(x)=\sum^{|\uk|+1}_{j=0}\big(\tau_j(x)\cdot g_j(x)\big)^{(\uk)}+\sum_{j>|\uk+1|}\dots$.
 Our bounds ensure that the infinite tail converges uniformly on the whole $\R^n$. Thus each $f^{(\uk)}$ is continuous.

\

\item
We claim: $\tau_j \cdot g_j -g_j\in \cc\cdot \cb_\infty$, for any $j$. Indeed, take $q\in \cb_\infty$,   satisfying $Z_\cb=q^{-1}(0)$.
 (See proposition \ref{Thm.Whitney.Thm.zeros.Strengthened}.)
As  $\tau_j \cdot g_j-g_j$ vanishes on $\cU_{\ep_j}(Z_\cb)$, the ratio $\frac{\tau_j \cdot g_j-g_j}{q}$ extends to a smooth function on $Ball_1(o)$.
Moreover, as $q$ is invertible on $Ball_1(o)\smin Z_\cb$, we get:  $\frac{\tau_j \cdot g_j-g_j}{q}\in \cc$.

 Therefore $\tau_j \cdot g_j-g_j\in  \ca\cdot \cb_\infty$.
Hence $f-\sum^N_{j=0} g_j\in I_N$, for any $N$.

 Thus the completion map sends $f$ to $\sum g_j$. Therefore the first row of \eqref{Eq.Completion.Surjective.Theorem} is exact.

\item (Exactness of the second row of \eqref{Eq.Completion.Surjective.Theorem})
Let $I_\bullet$ be a filtration of $C^\infty(\cU)$, as in the assumption. Take an element of the completion, $\sum g_j\in \widehat{C^\infty(\cU)}\ \!^{(I_\bullet)}$.
In the previous steps we have constructed a  representative $\tf\in C^\infty(\cU)$ of $\sum g_j$.

 We can assume $\tf>1$ on   $\cU\smin \cU(Z)$. Indeed, take some $\tau\in C^\infty(\cU)$ such that $\tau$ vanishes on a small neighbourhood of $Z$,
  but  $\tau>|\tf|+1$ on $\cU\smin \cU(Z)$.
  Then $\tau\in I_\infty$ and we replace $\tf$ by  $\tf+\tau$.

As the completion is surjective, $C^\infty(\cU)\twoheadrightarrow \widehat{C^\infty(\cU)}\ \!^{(I_\bullet)}$,
 the zero sets $V(I_j)$ stabilize, by lemma \ref{Thm.Zero.Loci.Stabilization}. Thus we assume
  $Z:=V(I_1)=V(I_2)=\cdots=V(I_\infty)$. By proposition \ref{Thm.Whitney.Thm.zeros.Strengthened}, $V(I_\infty)=V(\tau)$,
  for some  $\tau\in I_\infty$.

Take a sequence of positive reals, $\{\ep_i\to0\}$. Define
\beq
\cU_i:=\{x|\ \tau(x)^2>\ep_i\}\cap Ball_{\frac{1}{\ep_i}}(0)\cap \cU,\quad\quad
\text{ thus } \quad
 \cup \cU_i=\cU\smin Z.
\eeq

By the approximating  lemma \ref{Thm.Whitney.Approximation} there exists $f\in C^\infty(\cU)\cap C^\om(\cU\smin Z)$ satisfying: $|f-\tf|<\ep_i$ on $\cU_i$.
 Therefore $|f-\tf|<\tau^2$ on $\cU\smin Z$. Then $f-\tf\in I_\infty$. Therefore $f$ represents  $\sum g_j$, is analytic on $\cU\smin Z$
   and is positive on $\cU\smin \cU(Z)$.
\epr
\eee


\beR\label{Re.Broad.Class.of.Filtrations}
\bee[\bf i.]
\item It is peculiar that  the
sufficient condition of theorem \ref{Thm.Completion.Suff.Condition} seems to impose weaker restriction on
$I_\bullet$ than the necessary condition of theorem \ref{Thm.Completion.Nec.Condition}.
 The two theorems imply that in fact the two sets of conditions are equivalent.
\item
The class of filtrations allowed by the sufficient condition is
rather broad:
 \bei
 \item the condition is only on the germs of  $I_\bullet$,
 \item the ideals $\ca,\{\cc_k\},\{\cb_{k,j}\}$
 are not necessarily finitely/sub-analytically generated,
\eei
\item
Assuming only ``$I_j$ are polynomially generated" does  not
ensure the surjectivity, see example \ref{Ex.Ideal.Polynomially.Generated}.
\eee
\eeR

\bcor
 Assume $I_\bullet\sim I_\infty+I_\bullet\cap I(Z)^{\bl\bullet\br}$, for $Z=V(I_\infty)$. Then $\quots{R}{I_\infty}$ is complete for the filtration
 $\quots{I_\bullet}{I_\infty}$.
\ecor

\beR If the ideals $\{\cc_k\},\{\cb_{k,j}\}$ are analytically
generated then locally $\cc_k\cap \cb_k^{d_j}\sset \cc_k\cdot
\cb^{d_j-d_0}_k$, for some $d_0$. (By Artin-Rees in the Noetherian
ring $\R\{\ux\}$.) Thus the filtration of theorem
\ref{Thm.Completion.Suff.Condition}
 is equivalent to just $\{\ca+\sum_k \cc_k\cap \cb_k^\bullet\}$.

 For the (highly non-Noetherian) ring $\quots{C^\infty(\cU)}{J}$, with no analyticity assumption, no such simplification holds.
\eeR

\beR\label{Re.Positivity.Cannot.Strengthen}
For any neighbourhood $Z\sset \cU(Z)\sset \cU$ the theorem ensures a representative  $f\in C^\infty(\cU)\cap C^\om(\cU\smin Z)$
  that is positive on $\cU\smin \cU(Z)$.
 We cannot strengthen this to $f|_{\cU\smin Z}>0$, nor to $f$ is sign-semi-definite on $\cU(Z)$,
  due to the following immediate examples.
\bee[\bf i.]
\item
Let $I_\bullet=(x(x-1))^\bullet\sset C^\infty(\R^1)$, and
$g_1=x(x-1)(x-\frac{1}{2})$,
 $g_{\ge 2}=0$,
 $\hf=\sum g_j$.
  Suppose
  $f\in C^\infty(\R^1)$ represents $\hf$,
 then $f'|_0=\frac{1}{2}$ and $f'|_1=\frac{1}{2}$.
   Thus $f|_{(0,\ep)}>0$  and $f|_{(1-\ep,1)}<0$ for some $\ep>0$. Thus $f$ must have a zero on $(\ep,1-\ep)$.
   \item
   Let $I_\bullet=(y)^\bullet\sset C^\infty(\R^2)$ and $g_1=yx$, $g_{\ge2}=0$, $\hf=\sum g_i$. As before, for any $f$ representing $\hf$,
    and for any $x_0>0$, exists $\ep>0$
    such that $f>0$ on $\{x_0\}\times(0,\ep)$ and  $f<0$ on $\{-x_0\}\times(0,\ep)$. Thus $f$ necessarily vanishes at some point of $\cU(Z)\smin Z$,
     for any $Z\sset \cU(Z)$.
   \eee
\eeR

\section{Examples and applications}
\subsection{Examples} \label{Sec.Examples.surjectivity.of.completions}
  \bee[\bf i.]
\item As the simplest case suppose  $I_\bullet=\cm^\bullet\sset \quots{C^\infty(\R^n,o)}{J}, \quots{C^\infty(\cU)}{J}$, or, more generally,
 the filtration $I_\bullet$ is equivalent to $\cm^\bullet$.
 We get   the Borel lemma on the subscheme $V(J)\sseteq (\R^n,o),\cU$.
 The analyticity of the representative, $f\in C^w(\cU\smin\{0\})\cap C^\infty(\cU)$, recovers \cite[Theorem 1]{Shiota}.
\item
Suppose $V(I_j)=V(\cm)$ and $I_j\sseteq \cm^{d_j}$, with $d_j\to
\infty$, but $I_j\not\supseteq\cm^{N_j}$, for any $N_j<\infty$.
 (This happens, e.g. when $I_j$ is generated by flat functions.)
We still get the surjectivity of completion, though not
 implied by Borel lemma. (Now $I_\bullet$ is not equivalent to $\cm^\bullet$.)

For example, take a flat function $\tau\in \cm^\infty$, assume $V(\tau)=\{o\}\sset \R^n$. Then any formal   series  $\sum a_j\tau^j$, $a_j\in C^\infty(\cU)$,
 is presentable
 by some $f\in C^\infty(\R^n)$, in the sense: $f-\sum^N a_j\tau^j\in (\tau)^N$, for any $N\in \N$.

\item More generally, suppose $V(I_j)=V(I_1)=:Z$ and $I_j\sseteq  I(Z)^{\bl d_j\br}$, for a sequence $d_j\to\infty$. Again, theorem \ref{Thm.Completion.Suff.Condition}
 implies the surjectivity of completion. Note that we do not assume any regularity/subanalyticity conditions on the closed set $Z$.

If $Z\sset \cU$ is a discrete subset then we get a ``multi-Borel" lemma.

\item Take the ring $C^\infty(\R^n_x\times\R^m_y,o)$ with coordinates $x,y$, and the filtration $I_\bullet= (y)^\bullet$. The $(y)^\bullet$-completion map is the
Taylor map in $y$-coordinates, and theorem \ref{Thm.Completion.Suff.Condition} ensures its surjectivity:
\beq
C^\infty(\R^n\times\R^m,o)\twoheadrightarrow \widehat{C^\infty(\R^n\times\R^m,o)}\ ^{(y)^\bullet}=C^\infty(\R^n,o)[[y]].
\eeq
 Moreover, the preimage can be chosen $y$-analytic for $y\neq 0$.
This recovers the classical Borel  theorem, see   \cite[Theorem 1.2.6, pg. 16]{Hormander} and \cite[Theorem 1.3, pg. 18]{Moerdijk-Reyes}.

   \item
   For the ring $R=\quots{C^\infty\big((\R^n,o)\times[0,1]\big)}{J}$, we can interpret the elements as the families of function germs.
    Then we   get the surjectivity of completion in families,
  $R\twoheadrightarrow \hR$.
 For example, for $I_j=(\ux)^j$ we get a particular version of Borel lemma in families: any power series
$\sum a_{\um}(t) \ux ^\um$, with $ a_{\um}(t)\in C^\infty([0,1])$,
 is the $\ux $-Taylor expansion of some function germ $f_t( x )\in C^\infty\big((\R^n,o)\times[0,1]\big)$.

\item Many important filtrations are not by (differential) powers of ideals, neither are equivalent to this, see example \ref{Ex.Various.Filtrations}.
 As the simplest case, in Singularity Theory,  when studying the  germ (at the origin) of a non-isolated hypersurface singularity with singular
  locus  $\{x_1=0=x_2\}$, one often considers the filtration
\beq
I_\bullet=(x_1,x_2)^2\cdot(x_1,\dots,x_n)^\bullet\sset C^\infty(\R^n,o).
\eeq
More generally,
$I_j=x^2_1(x_1,y_1)^{d_{1,j}}+x^2_2(x_2,y_2)^{d_{2,j}}+\cdots +x^2_n(x_n,y_n)^{d_{n,j}}\sset C^\infty(\R^n_x\times \R^n_y,o)$
 is
  a typical filtration for complete intersections with non-isolated singularities.

This filtration satisfies the conditions of theorem \ref{Thm.Completion.Suff.Condition}. Thus one can use
  the surjectivity of completion  to pull-back various formal results  (over the completion)  to the $C^\infty$-statements.
 \eee

For further applications see \cite{Bel.Boi.Ker}, \cite{BGK}.

\subsection{An inverse Artin-Tougeron problem}\label{Sec.Inverse.Artin-Tougeron}
 Let $R=C^\infty(\cU_x\times \cU_y)$, for some open sets $\cU_x\sseteq \R^n_x$, $o\in \cU_y\sseteq \R^m_y$.
  Take the filtrations $I_\bullet$ of $C^\infty(\cU_x)$ and $(y)^\bullet$ of $C^\infty(\cU_y)$, and the corresponding filtration
   $\tI_\bullet:=R\cdot I_\bullet+R\cdot (y)^\bullet$. Identify $\hR^{(\tI_\bullet)}\isom{} \widehat{C^\infty(\cU_x)^{(I_\bullet)}}[[y]]$.
    An element $\hF\in (\hR^{(\tI_\bullet)})^q$ defines the system of (formal) power series equations, $\hF(y)=0$. Suppose this has a formal solution,
     $\hy(x)\in \big(\widehat{C^\infty(\cU_x)^{(I_\bullet)}}\big)^m$. The classical Artin-Tougeron problem asks for an ordinary solution,
 $\hF(y(x))=0$, for some $y(x)\in \big(C^\infty(\cU_x)\big)^m$ that is sent by the completion to $\hy(x)$.

  The inverse Artin-Tougeron problem (see e.g. \cite{Hubl} and references therein) asks for $F\in R^q$ that is sent to $\hF$,
   and $y(x)\in (C^\infty(\cU_x))^m$ that is sent to $\hy$, such that $F(y(x))=0$.

In the $C^\infty$-case this is resolved now trivially as follows. Assume the completion maps are surjective,
 $C^\infty(\cU_x)\twoheadrightarrow\widehat{C^\infty(\cU_x)^{(I_\bullet)}}$ and $R\twoheadrightarrow \hR$.
   Take any representatives, $\tF$ of $\hF$ and $y(x)$ of $\hy(x)$. Then $\tF(y(x))\in \tI_\infty\cdot (C^\infty(\cU_x))^q$. Therefore
    $F(y):=\tF(y)-\tF(y(x))$ is the needed representative of $\hF$, and $F(y(x))=0$.

\subsection{$C^\infty$-contraction of ideals}\label{Sec.Contraction.of.Ideals}
 Given an ideal $\ca\sset R$ and a ring homomorphism $R\stackrel{\phi}{\to}S$, one often takes the contraction,
 $\phi^{-1}(S\cdot\phi(\ca))\sset R$. Then $\phi^{-1}(S\cdot\phi(\ca))\supseteq \ca+Ker(\phi)$, the inclusion can be proper.
 The equality holds for the completion of Noetherian local rings  \cite[Chapter III]{Bourbaki}, but fails otherwise.
 In our case the surjectivity of completion forces the equality.
 \bcor
If $R\stackrel{\phi}{\twoheadrightarrow}\hR^{(I_\bullet)}$ then $\phi^{-1}\big(\hR^{(I_\bullet)} \phi(\ca)\big)= \ca+Ker(\phi)$, for any ideal $\ca\sset R$.
 \ecor
 Equivalently, $\phi(\ca)=\hR^{(I_\bullet)}\phi(\ca)$. Thus any ideal in $\hR^{(I_\bullet)}$ is of the form $\phi(\ca)$, with $V(\ca)\sseteq Z$.
 This answers a question of N. Zobin, \cite{Fefferman}, initially asked for the filtrations $I(Z)^{\bl \bullet\br}$.
\bpr
For any elements $\{a_i\}$ of $\ca$ and $\{\hf_i\}$ of $\hR^{(I_\bullet)}$ the preimage of $\sum \hf_i\cdot\phi(a_i)$ is  $\sum f_i\cdot a_i$ .
\epr

\subsection{Surjectivity for completion of modules}\label{Sec.Surjectivity.Completion.of.Modules}
Let $M$ be a module over the ring $R=\quots{C^\infty(\cU)}{J}$. A filtration $I_\bullet$ of $R$ induces the filtration $M_\bullet=I_\bullet\cdot M$.
\bcor
\bee[1. ]
\item
If $M$ is finitely generated and $I_\bullet$  satisfies the conditions of theorem \ref{Thm.Completion.Suff.Condition}
 then the completion map $M\to \hM^{(I_\bullet)}$ is surjective.
\item
 Assume $M$ contains an $R$-regular element, i.e. $Ann_R(z)=0$ for some $z\in M$. If the map $M\to \hM^{(I_\bullet)}$ is surjective then
  the conditions of theorem \ref{Thm.Completion.Nec.Condition} are satisfied.
\eee
\ecor
\bpr
{\bf 1.} Fix some generators $\{z_\al\}$ of $M$, then an element of $\hM^{(I_\bullet)}$ is presentable as $\sum_\al (\sum_j g_{j,\al} z_\al)$.
 Here $g_{\bullet,\al}\in I_\bullet$. Thus if $f_\al\in R$ is a representative of $\sum g_{\bullet,\al}$ then $\sum_\al f_\al z_\al$ represents this element.

{\bf 2.} Let $z\in M$ be a regular element then $R\cong R\cdot\{z\}\sset M$. Thus the surjectivity of $M\to \hM^{(I_\bullet)}$ implies the surjectivity of
 $R\to \hR^{(I_\bullet)}$.
\epr
\beR The finitely-generated/regular assumptions of this corollary are necessary.
\bee[\bf i.]
\item Let $R=C^\infty(\cU)$ with the filtration by powers of a maximal ideal, $\cm^\bullet$.
 Take an infinite direct sum, $M=\oplus R\bl t_i\br$,  with generators $\{t_i\}$.
  Take an element $0\neq x\in C^\infty(\cU)$. Then $\hM^{(\cm)}\ni \sum^\infty_{n=0}x^n \Big(\sum^n_{j=1} t_j\Big)$. This element has no $M$-representative, as $M$
  contains only finite sums, $\sum^k_{j=1} f_j t_j$, for $f_j\in C^\infty(\cU)$.
  \item Let $R=C^\infty(\R^n)$, $M=\quots{R}{(f)}$, for some $0\neq f\in R$. Take a filtration  $I_\bullet$  satisfying: $I_N\sseteq (f)$ for $N\gg1$.
   Then $\hM^{(I_\bullet)}\cong M$, as $R$ and $\hR^{(I_\bullet)}$-modules, but the completion  $R\to \hR^{(I_\bullet)}$ is not necessarily surjective.
\eee
\eeR

\subsection{Lifting $\hR^{(I_\bullet)}$-modules to $R$-modules}\label{Sec.Lifting.of.Modules}
 Let $R=C^\infty(\cU)$, with a filtration $I_\bullet$.
 Denote $Z:=V(I_\infty)$ and define the $Z$-support of an  $\hR^{(I_\bullet)}$-module $\hM$ by $Supp_Z(\hM):=\{x\in Z|\ M|_{(\cU,x)}\neq 0\}$.
\bcor
Assume $I_\bullet$ satisfies the conditions of theorem \ref{Thm.Completion.Suff.Condition}. Suppose an  $\hR^{(I_\bullet)}$-module $\hM$ is finitely-generated
 (resp.  finitely-presented). Then there exists a finitely-generated  (resp.  finitely-presented) $R$-module $M$ satisfying:
  $M\otimes \hR^{(I_\bullet)}\cong\hM$ and $Supp(M)=Supp_Z(\hM)$.
\ecor
This extends (partially) \cite{Tougeron1976-2}.

\bpr
Choose a presentation $\hR^n\stackrel{\hA}{\to}\hR^m\to \hM\to0$. Let $\tA\in \Mat$ be an $R$-representative of $\hA$,
 then $Coker(\tA)\otimes\hR\approx Coker(\tA\otimes\hR)=\hM$.
 To ensure the support take some $\tau\in I_\infty$ satisfying $Z=V(\tau)$,
  see proposition \ref{Thm.Whitney.Thm.zeros.Strengthened}. Define $A=\bbm \tA|\tau\one\ebm\in Mat_{m\times(m+n)}(R)$ and
   $ R^{n+m}\stackrel{ A}{\to} R^m\to  M\to0$. Then $M\otimes\hR\cong \hM$, and $Supp(M)\sseteq Z$. Finally, for each point $x\in Z$ the map $A|_x$
    is surjective iff $\hA|_x$ is surjective. Thus  $Supp(M)=Supp_Z(\hM)$.
\epr

\section{$C^\infty$-representatives with further properties}

\subsection{Representatives with prescribed decay on $\di\cU$}\label{Sec.Surjectivity.Prescribed.Decay}
For various applications one wants a representative with fast decay ``at infinity".
 Fix a (multi-)sequence of functions $q_\uk\in C^\infty(\cU)$ satisfying:  $\{q_\uk> 0\}_\uk$ on $\cU$, and  $\{q_\uk\to 0\}_\uk$  for $x\to \di\cU$.
Accordingly  define the Schwartz-type vector space of functions decaying with their derivatives,
 \beq
 Sch_{\di\cU}:=\big\{f\in C^\infty(\cU)|\quad f^{(\uk)}=O(q_\uk),\ \text{as }x\to x_0,\ \forall x_0\in \di\cU,\ \forall \uk\big\}.
 \eeq
 For any  $\hf\in \widehat{C^\infty(\cU)}^{(I_\bullet)}$
we want to ensure a representative   with the prescribed decay on $\di\cU$, i.e. to
strengthen the surjectivity of completion to $ Sch_{\di\cU}\twoheadrightarrow \widehat{C^\infty(\cU)}^{(I_\bullet)}$.

 This cannot hold without further assumptions, e.g. $Sch_{\di\cU}$ is not a unital ring, it does not contain $1$.
 More generally, let $Z=V(I_\infty)$ and take the element $g_1+\sum 0\in \widehat{C^\infty(\cU)}^{(I_\bullet)}$, where
    $\liml_{x\to x_0}g^{(\uk)}_1(x)\neq0$,
 for some $x_0\in \overline{Z}\cap \di\cU$ and some $\uk$.  (Here $\overline{Z}$ is the closure in $\R^n$.)
  All the $C^\infty(\cU)$-representatives of $g_1$ are of the form $\{g_1\}+I_\infty$.
 None of these belongs to $ Sch_{\di\cU}$.

Therefore we should a priori impose the needed behaviour on the elements of $C^\infty(\cU)$, at least at the points of $\overline{Z}\cap \di\cU$.
 Accordingly, we define the  space of functions whose germs at the points of $\overline{Z}\cap \di\cU$ have the prescribed decay:
\beq
Sch_{\overline{Z}\cap \di\cU}:=\big\{f\in C^\infty(\cU)|\  \text{ for any } x\in \overline{Z}\cap \di\cU
 \text{ the germs at $x$ satisfy:}  \ [f]_x\in  [Sch_{\di\cU}]_x
  \big\}.
\eeq
Thus $Sch_{\overline{Z}\cap \di\cU}\supseteq Sch_{\di\cU}$.
If $\overline{Z}\cap \di\cU=\empty$ then $Sch_{\overline{Z}\cap \di\cU}=C^\infty({\cU})$.
 If $Z=\cU$ then $Sch_{\overline{Z}\cap \di\cU}=Sch_{\di\cU}$.

\

Neither   $Sch_{\overline{Z}\cap \di\cU}$ nor $Sch_{\di\cU}$ are  unital rings. Still, for a filtration $I_\bullet\sset Sch_{\overline{Z}\cap \di\cU}$ with
 $V(I_\infty)=Z$, we can take the completion of vector spaces,
 \beq
\widehat{Sch_{\overline{Z}\cap \di\cU}}^{(I_\bullet)}:=\liml_{\leftarrow}\quot{Sch_{\overline{Z}\cap \di\cU}}{I_\bullet}.
 \eeq
We establish the surjectivity   with representatives of prescribed decay rate on $\di\cU$.
\bprop
Take a  filtration $I_\bullet\sset Sch_{\overline{Z}\cap \di\cU}$, with $V(I_\infty)=Z$. Suppose $I_\bullet$ satisfies the
 assumptions of theorem \ref{Thm.Completion.Suff.Condition}. Then the following rows are exact:
 \beq
\bM 0\to I_\infty\to&Sch_{\overline{Z}\cap \di\cU} &\stackrel{}{\to}& \widehat{Sch_{\overline{Z}\cap \di\cU}}^{(I_\bullet)} &\to 0
\\& \cup
&& || \quad
\\
&  Sch_{\di\cU}\cap C^\om(\cU\smin Z)& \to &\widehat{Sch_{\overline{Z}\cap \di\cU}}^{(I_\bullet)}&\to 0
\eM
\eeq
Moreover, for any neighbourhood $Z\sset \cU(Z)$ and any element $\hf\in \widehat{Sch_{\overline{Z}\cap \di\cU}}^{(I_\bullet)}$ there exists a representative
 $f\in   Sch_{\di\cU}\cap C^\om(\cU\smin Z)$ that satisfies  $f>0$ on $\cU\smin \cU(Z)$.
\eprop
\bpr
 We can assume that $\cU\sset \R^n$ is bounded. (Recall, any open subset of $\R^n$ is real-analytically equivalent to a bounded open,
  e.g. by $x\to \frac{x}{1+||x||^2}$.)

Take an element $\sum g_j\in \widehat{Sch_{\overline{Z}\cap \di\cU}}^{(I_\bullet)}$.
 Here $g_j$ does not necessarily belong to $Sch_{\di\cU}$. However,
  at each point $x\in \overline{Z}\cap \di\cU$
  the inclusion holds for germs, $[g_j]_x\in [Sch_{\di\cU}]_x$.
 Therefore for each $j$ we take a small (relatively) open neighbourhood
 $\overline{Z}\sset \cU_j(\overline{Z})\sset \overline{\cU}$
 such that $g_j|_{\cU_j(\overline{Z})}\in I_j\cap Sch_{\di\cU}|_{\cU_j(\overline{Z})}$.

  Take an element $\tau_j\in I_\infty$ that satisfies     $\tau_j|_{\overline{\cU}\smin \cU_j(\overline{Z})}=1$.
   It exists e.g. by proposition \ref{Thm.Whitney.Thm.zeros.Strengthened}.

   Then the element  $\sum g_j\in \widehat{Sch_{\overline{Z}\cap \di\cU}}^{(I_\bullet)}$ is equivalent to
    $\sum g_j(1-\tau_j)$. 
   By construction:
   $g_j(1-\tau_j)\in I_j\cap Sch_{\di\cU}$, in particular this is flat on $\di\cU$.

   Extend    $g_j(1-\tau_j)$ to $\tg_j\in C^\infty(\R^n)$ by zero outside of $\cU$.
    Take the filtration $I_\bullet\cap\cq\cap I(\R^n\smin \cU) $ of the ring $C^\infty(\R^n)$,
    here $\cq$ is the ideal defined by the $\{q_\uk\}$-vanishing conditions.
 This filtration satisfies the assumptions of theorem  \ref{Thm.Completion.Suff.Condition}. Thus we have a representative
  $\tf\in C^\infty(\R^n)\cap C^w(\cU\smin Z)$ of $\sum\tg_j$. Moreover, for any given neighbourhood $Z\sset\cU(Z)\sset \cU $, we can assume $\tf>0$ on $\cU\smin \cU(Z) $.
Finally,  the restriction $f:=\tf|_\cU\in Sch_{\di\cU}\cap C^w(\cU\smin Z)$ is the needed representative of  $\sum g_j\in \widehat{Sch_{\overline{Z}\cap \di\cU}}^{(I_\bullet)}$.
\epr
\bex
Assume $Z\sset \cU$ is compactly embedded, then $Sch_{\overline{Z}\cap\di\cU}\!=\!C^\infty(\cU)\!=:\!R$. We get: any $\hf\!\in\! \hR^{(I_\bullet)}$\!
 admits a representative $f\!\in\! Sch_{\di\cU}\!\cap\! C^w(\cU\!\smin\! Z)$,  decaying (at the prescribed rate) on $\di\cU$.
\eex

\subsection{Representatives with imposed $\R$-linear conditions}\label{Sec.Surj.with.Linear.Conditions}
For various applications one needs $C^\infty$-representatives of  power series with prescribed integrals, e.g. $\{\int^\infty_{-\infty} f\cdot x^k dx=c_k\}$
 for the Hamburger moment problem. More generally, take    a (possibly infinite)  set $\{l_k\}$ of $\R$-linear functionals defined on
 $R=C^\infty(\cU)$, or on a vector subspace of $R$.
 (We do not assume that $\{l_k\}$ are continuous in whichever sense.)
\bed
1. The conditions $\{(l_k,c_k)\}$ are resolvable   if exists $f\in R$ satisfying $\{l_k(f)=c_k\}$.

2. The conditions $\{(l_k,c_k)\}$ are $\hf$-resolvable   if exists a representative $f\in R$ of $\hf$ that satisfies $\{l_k(f)=c_k\}$.
\eed

Resolvability of $\{(l_k,c_k)\}$ can be obstructed in various ways.
 \bex
\bee[\bf i.]
\item
    Define $g_\infty(x):=\sum g_k(x)\in C^0(\cU)$, suppose the series converges uniformly  on $\cU$.
     Define $l_k(f):=\int_\cU g_k fd\cU$, for $k\le \infty$, these functionals act on the subspace  of functions for which
       these integrals converge. Then, $\sum l_k=l_\infty$,
     imposing restrictions on $\{c_k\}$.

\item
  Take a convergent sequence of points $x_k\to x_\infty\in \cU$, assume $\{x_k\}_k$ are pairwise distinct for $k\le \infty$.
Define $l_k(f):=f(x_k)$ for $k\le \infty$.
Thus $l_k\to l_\infty$, pointwise on $R$.
 The functionals $\{l_k\}$ are   linearly independent. But the conditions $\{(l_k,c_k)\}$
   are resolvable iff $\lim c_k=c_\infty$. 
\item
 Take a sequence of points $\{x_k\}$  in $\cU$. Define
 $l_k(f):=\int_\cU f\cdot g_kdx$, where $\{g_k\}$
 are bump functions on $\{Ball_{\frac{1}{k!}}(x_k)\}$.
 Then $\{l_k\}$ are linearly independent.
  Assume the sequence $\{x_k\}$ has a condensation point $x_0\in \cU$, then the conditions $\{l_k(f)=1\}$ force $f$ to
  explode at $x_0$.
\item Assume $\cU$ is bounded and take an exhaustion $\cU_1\sset \cU_2\sset\cdots\sset\cup\cU_k=\cU$,
 by compactly embedded sets. Define
 $l_k(f):=\int_\cU f\cdot (1-g_k)dx$, where $\{g_k\}$
 are bump functions on $\{\cU_k\}$. Assume $\cU\smin \cU_k$ decreases fast. Then the conditions $\{l_k(f)=1\}$ force $f$
   to explode at some point(s) of $\di\cU$.
   While $f$ remains $C^\infty$, the conditions $l_k(f)=c_k$  become ill-defined, as the integrals diverge.
\eee
 \eex
Therefore any statement on  $\hf$-resolvability of $\{(l_k,c_k)\}$ must assume the resolvability of $\{(l_k,c_k)\}$.
\bel
The conditions $\{(l_k,c_k)\}$ are $\hf$-resolvable iff   $\{(l_k,c_k)\}$ are resolvable and moreover $\cap Ker(l_k)\twoheadrightarrow \hR$.
\eel
\bpr $\Rrightarrow$ Take $\{c_k=0\}$ and any $\hf\in \hR$.

$\Lleftarrow$ Take some $g\in R$ satisfying $\{l_k(g)=c_k\}$. Take some $\tf\in \cap Ker(l_k)$ whose image in $\hR$ is $\hat{\tf}=\hf-\hg$.
 Then $g+\tf$ is the needed solution.
\epr
\bex
 \bee[\bf i.]
  \item
Assume the completion is surjective, $R\twoheadrightarrow \hR^{(I_\bullet)}$, denote $Z=V(I_\infty)$.
 Call a functional ``punctual" if it is presentable as $l(f):=\sum_j c_{\um_j}\cdot f^{(\um_j)}|_{x_j}$, the values of prescribed partials at prescribed points.
   If this sum is infinite then $l$ acts on a subspace of $R$ for which the sum converges.
  Call the points $\{x_j\}$ ``the base points of $l$".

  Assume the collection  $\{l_k\}$ consists of
 an arbitrary number of punctual functionals, all of whose base-points lie in $\cU\smin \cU(Z)$, for a neighbourhood $\cU(Z)$.
 Then $\cap Ker(l_k)\twoheadrightarrow \hR$. Hence $\{(l_k,c_k)\}$ are $\hf$-resolvable iff $\{(l_k,c_k)\}$ are resolvable.

  \item   Assume the collection  $\{l_k\}$ consists of
  an arbitrary number of punctual functionals, with base-points   in $\cU\smin \cU(Z)$,   and of
  a finite number of moments, $l_k(f)=\int_\cU h_k(x)\cdot f(x)d^nx$, for some  functions $\{h_k\}$ integrable on $\cU$.
Then $\cap Ker(l_k)\twoheadrightarrow \hR$. Hence $\{(l_k,c_k)\}$ are $\hf$-resolvable iff $\{(l_k,c_k)\}$ are resolvable.
  \eee

\eex
\beR
One cannot impose an infinite collection of moments without non-trivial assumptions on $\{c_k\}$. Even if the conditions
  $\int_\cU h_k\cdot f\cdot d^nx=c_k$
 are resolvable on $R$, the solution can be unique. (This can happen e.g. for Hamburger moment problem.) Then no $\hf$-solvability
   of $\{(l_k,c_k)\}$ can be achieved.
\eeR

\end{document}